# The solution infinite horizon noncooperative differential game with nonlinear dinamics without the Hamilton-Jacobi-Bellman equation.

Jaykov Foukzon

Israel Institute of Technology

jaykovfoukzon@list.ru

**Abstract:** For a non-cooperative m-persons differential game, the value functions ofthe various players satisfy a system of Hamilton-Jacobi-Bellman equations.Nashequilibrium solutions in feedback form can be obtained by studying a related system of P.D.E's.A new approach, which is proposed in this paper allows one to construct the feedback optimal control $\alpha^*(x) = (\alpha_1^*(x),\ldots,\alpha_m^*(x))$ and cost functions $\mathbf{J}_i(t,x_0), i = 1,\ldots,m$ directly,i.e.,without any reference to the corresponding Hamilton-Jacobi-Bellman equations.

# I.Introduction.

A dynamic game is a system with the following attributes:

**(a)** It has $m$ persons, players, or decision-makers.
**(b)** Player $i$ chooses a control variable $\mathbf{u}_i$ from a set of admissible controls $\mathbf{U}_i$.
**(c)** It has a time horizon which is defined by the interval $[t_0, t_f]$ where $t_0$ is known and fixed, and $t_f$ may be fixed or free and it may be finite: $t_f < \infty$ or infinite: $t_f \leq \infty$.
**(d)** It has a state $x(t)$ at time $t$, $t \in [t_0, t_f]$ which is an element of a finite dimensional vector space $\mathbb{R}^n$. The evolution of the state is such that $x(t)$ is uniquely determined by the values of $\mathbf{u}_i$ on $[t_1, t]$, $i = 1,\ldots,m$ and $x(t_1)$ for any $t_1$, satisfying $t_0 < t_1 < t$.

We only consider state evolutions describable by differential equations.

**(e)** Each player $i = 1,\ldots,m$ has a real scalar cost function $J_i$ which is a mapping from $\mathbb{R}^n$ and $\mathbf{U}_i, i = 1,\ldots,m$ to the set of real numbers $\mathbb{R}$.

**(f)** Each player $i = 1,\ldots,m$ has knowledge of an information set $\mathbf{I}_i(t)$ which may include the differential equations for state evolution, the state v, its own cost function mapping as well as those of the other players, and control strategies of the other players. The set $\Im = \{\mathbf{I}_i(t)|i = 1,\ldots,m, t \in [t_0, t_f]\}$ is called the information structure of the game.

**(g)** Each player $i$ has a control law or strategy $\breve{u}_i(t) : \mathbf{I}_i(t) \to \mathbf{U}_i$ which is a mapping from the information set $\mathbf{I}_i(t)$ to the control space $\mathbf{U}_i$.
A dynamic game whose state evolution is given by a differential equation is called a differential game.

The problem we wish to consider is an $m$-player dynamic game in which the state of the $j$-th player at time $t \in [t_0, t_f]$, is described by a dynamical system of the form:

$$\dot{x}_i(t) = f_i(t, \mathbf{x}(t), \mathbf{u}(t)),$$

$$x(t_0) = x_0, x(t_f) = x_f, \tag{1.1}$$

where $\mathbf{x}(t) = (x_1(t), x_2(t),\ldots,x_k(t)) \in \mathbb{R}^{n_1} \times \mathbb{R}^{n_2} \times \ldots \times \mathbb{R}^{n_k} = \mathbb{R}^N$, i.e. for each $i = 1,\ldots,m$, $x_i(t) \in \mathbb{R}^{n_i}$ and $n_1 + n_2 + \ldots + n_k = N$. Also, we let

$$\mathbf{u}(t) = (u_1(t), u_2(t),\ldots,u_m(t)).$$

We assume that these strategies are defined as functions $u_i(t) : [t_0, t_f] \to \mathbb{R}^{n_i}$ and satisfy control constraints of the form:

$$u_i(t) \in \mathbf{U}_i(t) \subsetneq \mathbb{R}^{n_i} \quad a.e. \ t_0 \le t \le t_f, \tag{1.2}$$

where for each $t \in [t_0, t_f]$ and each $i$ the set $\mathbf{U}_i(t)$ is assumed to be closed and nonempty. Finally we assume that the sets:

$$M_i \triangleq \{(t, \mathbf{x}, u_i) \in [t_0, t_f] \times \mathbb{R}^N \times \mathbb{R}^{n_i} | u_i(t) \in \mathbf{U}_i(t)\} \tag{1.3}$$

are closed and nonempty for each $i = 1, 2,\ldots,m$ and that $f_i(\circ,\circ,\circ) : M_i \to \mathbb{R}^{n_i}$ are continuous.

The performance of each player is measured through a performance index described by an integral of the form:

$$\mathbf{J}_i(\mathbf{x}(\circ), u_i(\circ)) = \int_{t_0}^{t_f} f_i(t, \mathbf{x}(t), u_i(t)) dt \qquad (1.4)$$

It is assumed that $\mathbf{I}_i$ for each $i$ includes knowledge of $f = (f_1, f_2, \ldots, f_m)$.

This paper deal with the case that there is more than one player. For ease of exposition we will just deal with the two-player case.

# II.1. Infinite Horizon 2-Persons noncooperative dissipative differential game with nonlinear dynamics.

In this section we consider a scalar 2-persons differential game, with nonlinear dynamics

$$\dot{x} = f_1(x, \alpha_1) + f_2(x, \alpha_2(x)),$$

$$x(0) = y \in \mathbb{R}, \qquad (2.1.1)$$

$$\alpha_i(t) \in A_i.$$

The functions $t \mapsto \alpha_i(t)$, $i = 1, 2$, represent the controls implemented by the $i$-th player, chosen within a compact set of admissible controls $\mathcal{A}_i \subset \mathbb{R}$. The game takes place on $[0, +\infty)$ and each player is subject to a running cost, exponentially discounted, of the following form:

$$\mathbf{J}_i(x(t), \alpha_i(t)) \triangleq \int_0^\infty e^{-t} \left[ h_i(x(t)) + \frac{\alpha_i^2(t)}{2} \right] dt. \qquad (2.1.2)$$

Assume here that both $h_i(x)$ are piecewise smooth functions with bounded derivatives.

A couple of feedback strategies $(\alpha_1^*(x), \alpha_2^*(x))$ represents a *Nash equilibrium solution* for the game (2.2.1)–(2.2.2) if the following holds. For $i \in \{1,2\}$, the feedback control $\alpha_i = \alpha_i^*(x)$ provides a solution to the the optimal control problem for the $i$-th player,

$$\min_{\alpha_i(\cdot)} \mathbf{J}_i(\alpha_i), \qquad (2.1.3)$$

where the dynamics of the system is

$$\dot{x} = f_1(x, \alpha_1^*) + f_2(x, \alpha_2^*(x)),$$

$$x(0) = y \in \mathbb{R}, \qquad (2.1.4)$$

$$\alpha_i(t) \in A_i, \ i = 1, 2.$$

More precisely, we require that, for every initial data $y \in \mathbb{R}$, the Cauchy problem

$$\dot{x} = f_1(x, \alpha_1(x)) + f_2(x, \alpha_2^*(x)),$$

$$x(0) = y, \qquad (2.1.5)$$

should have at least one Caratheodory solution $t \mapsto x(t)$, defined for all $t \in [0, \infty)$.
Moreover, for every such solution and each $i = 1, \ldots, m$, the cost to the $i$-th player should provide the minimum for the optimal control problem (2.1.4)-(2.1.5). We recall that a Caratheodory solution is an absolutely continuous function $t \mapsto x(t)$ which satisfies the differential equation (2.1.5) at almost every $t > 0$.

## II.2. Infinite Horizon 2-Persons noncooperative differential game with linear dynamics.

In this section we consider a scalar 2-persons differential game, with linear

dynamics:

$$\dot{x} = \alpha_1 + \alpha_2,$$

$$x(0) = y \in \mathbb{R}. \tag{2.2.1}$$

$$\alpha_i \in \mathbb{R}.$$

The functions $t \mapsto \alpha_i(t)$, $i = 1, 2$, represent the controls implemented by the $i$-th player, chosen within a compact set of admissible controls $\mathcal{A}_i \subset \mathbb{R}$. The game takes place on $[0, +\infty)$ and each player is subject to a running cost, exponentially discounted, of the following form:

$$\mathbf{J}_i(x(t), \alpha_i(t)) \triangleq \int_0^\infty e^{-t} \left[ h_i(x(t)) + \frac{\alpha_i^2(t)}{2} \right] dt. \tag{2.2.2}$$

Assume here that both $h_i$ are piecewise smooth functions with bounded derivatives.

A couple of feedback strategies $(\alpha_1^*(x), \alpha_2^*(x))$ represents a *Nash equilibrium solution* for the game (2.2.1)–(2.2.2) if the following holds. For $i \in \{1, 2\}$, the feedback control $\alpha_i = \alpha_i^*(x)$ provides a solution to the the optimal control problem for the $i$-th player,

$$\min_{\alpha_i(\cdot)} \mathbf{J}_i(\alpha_i), \tag{2.2.3}$$

where the dynamics of the system is

$$\dot{x} = \alpha_i^* + \alpha_j^*(x),$$

$$\alpha_i(t) \in A_i, j \neq i. \tag{2.2.4}$$

More precisely, we require that, for every initial data $y \in \mathbb{R}$, the Cauchy problem

$$\dot{x} = \alpha_1^*(x) + \alpha_2^*(x),$$

$$x(0) = y,$$

(2.2.5)

should have at least one Caratheodory solution $t \mapsto x(t)$, defined for all $t \in [0, \infty)$.

Moreover, for every such solution and each $i = 1, \ldots, m$, the cost to the $i$-th player should provide the minimum for the optimal control problem (2.2.1)-(2.2.1). We recall that a Caratheodory solution is an absolutely continuous function $t \mapsto x(t)$ which satisfies the differential equation in (2.2.5) at almost every $t > 0$.

The vector function $u(x) = (u_1(x), u_2(x))$ thus satisfies the stationary system of equations:

$$u_i(x) = H_i(x, u_1', u_2'),\qquad(2.2.6)$$

where the Hamiltonian functions $H_i, i \in \{1, 2\}$ are defined as follows. For each $p_j \in \mathbb{R}$, assume that there exists an optimal control value $\alpha_j^*(x, p_j)$ such that

$$p_j \cdot \alpha_j^*(x, p_j) + \psi_j(x, \alpha_j^*(x, p_j)) = \min_{a \in A_j} \{p_j \cdot a + \psi_j(x, a)\},$$

$$j \in \{1, 2\}.$$

(2.2.7)

Then

$$H_i(x, p_1, p_2) \triangleq p_i \cdot \alpha_j^*(x, p_j) + \psi_i(x, \alpha_i^*(x, p_i)).\qquad(2.2.8)$$

for $i, j \in \{1, 2\}$ and $i \neq j$.

In general, even in cases as easy as $\psi_i = \alpha_i^2/2$, this system will have infinitely many solutions defined on the whole $\mathbb{R}$. And not every solution corresponds to a Nash equilibrium for the initial game. To single out a (hopefully unique) admissible solution, and therefore a Nash equilibrium for the differential game, additional requirements must be imposed [16]:

**Definition**  2.2.1. Namely a solution $u$ to (2.2.6) is said to be an admissible solution if the following holds [16]:

(A1)  $u(x)$ is absolutely continuous and its derivative $u'(x)$ satisfies (2.2.6) at a.e. point $x \in \mathbb{R}$.

(A2)  $u(x)$ has sublinear growth at infinity; namely, there exists a constant $C > 0$ such that, for all $x \in \mathbb{R}$,

$$|u(x)| \leq C(1 + |x|). \tag{2.2.9}$$

(A3)  At every point $y \in R$, the derivative $u'$ admits right and left limits $u'(y+), u'(y-)$ and at points where $u'$ is discontinuous, these limits satisfy at least one of the conditions:

$$u'_1(y+) + u'_2(y+) \leq 0$$

or  (2.2.10)

$$u'_1(y-) + u'_2(y-) \geq 0.$$

Because of the assumption on $h'_i(x)$, the cost functions $h_i(x)$ are Lipschitz continuous. It is thus natural to require the value functions $u_i$ to be absolutely continuous, with sub-linear growth as $x \to \pm\infty$. The motivation for the assumption (A3) is quite simple. Observing that, in (2.2.5), the feedback controls are $\alpha^*_i = -u'_i(x)$, the condition (2.2.10) provides the existence of a local solution to the Cauchy problem

$$\dot{x} = -u'_1(x) - u'_2(x), x(0) = y \tag{2.2.11}$$

forward in time. In the opposite case, solutions of the O.D.E. would approach $y$ from both sides, and be trapped. Thus in general, system (2.2.6) will have infinitely many solutions. To single out a (hopefully unique) admissible solution, corresponding to a Nash equilibrium for the differential game, additional requirements must be imposed.
These are of two types:

 (i) Asymptotic growth conditions as $|x| \to \infty$.

(ii) Jump conditions, at points where the derivative $u'(x)$ is discontinuous.

The general system of H-J equations (2.2.6) for the value functions now takes the form

$$u_1(x) = h_1(x) - u'_1 u'_2 - (u'_1)^2/2,$$

$$u_2(x) = h_2(x) - u'_1 u'_2 - (u'_2)^2/2. \tag{2.2.12}$$

and the optimal feedback controls are given by

$$\alpha^*_i(x) = -u'_i(x). \tag{2.2.13}$$

Differentiating (2.2.12) and setting $p_i = u'_i$ one obtains the system

$$h'_1 - p_1 = (p_1 + p_2)p'_1 + p_1 p'_2,$$

$$h'_2 - p_2 = p_2 p'_1 + (p_1 + p_2)p'_2, \tag{2.2.14}$$

Set

$$\Lambda(p) \triangleq \begin{pmatrix} p_1 + p_2 & p_1 \\ p_2 & p_1 + p_2 \end{pmatrix},$$

$$\Delta(p) \triangleq \det \Lambda(p), \tag{2.2.15}$$

From (2.2.14)-(2.2.15) we deduce

$$p'_1 = [\Delta(p)]^{-1}[-p_1^2 + (h'_1 - h'_2)p_1 + h'_1 p_2],$$

$$p'_2 = [\Delta(p)]^{-1}[-p_2^2 + (h'_2 - h'_1)p_2 + h'_2 p_1].$$

(2.2.16)

Notice that

$$\frac{1}{2}(p_1^2 + p_2^2) \leq \Delta(p) \leq 2(p_1^2 + p_2^2).$$

(2.2.17)

In particular, $\Delta(p) > 0$ for all $p = (p_1, p_2) \neq (0, 0)$. Hence, $\Lambda(p)$ is invertible outside the origin and, for $p \neq (0, 0)$, we can restrict the study to the equivalent system

$$p'_1 = (h'_1 - h'_2)p_1 + h'_1 p_2 - p_1^2,$$

$$p'_2 = (h'_2 - h'_1)p_2 + h'_2 p_1 - p_2^2,$$

(2.2.18)

For piecewise smooth solutions, jumps are only allowed from any point

$(p_1^-, p_2^-)$ with

$$0 \leq p_1^- + p_2^-$$

(2.2.19)

to the symmetric point

$$(p_1^+, p_2^+) = (-p_1^-, -p_2^-).$$

(2.2.20)

Consider the game for two players, with dynamics

$$\dot{x} = \alpha_1 + \alpha_2, x(0) = y. \quad (2.2.21)$$

and cost functionals

$$J_i(\alpha_i) = \frac{1}{2} \int_0^\infty e^{-t} \alpha_i^2(t) dt, \quad (2.2.22)$$

$$i = 1, 2.$$

The system of H-J (2.12) takes the simple form

$$u_1 = -\frac{1}{2}(u_1')^2 - u_1' u_2',$$
$$u_2 = -\frac{1}{2}(u_2')^2 - u_1' u_2' \quad (2.2.23)$$

The obvious admissible solution is $u_1 \equiv u_2 \equiv 0$, corresponding to identically zero controls, and zero cost. We now observe that the functions

$$u_1(x) = -\frac{1}{2}x^2, u_2(x) = 0 \quad (2.2.24)$$

provide solution, which does not satisfy the growth conditions (2.2.9). In these case, the corresponding feedback:

$$\alpha_i^*(x) = -u_i(x), i = 1, 2,$$
$$\alpha_1^*(x) = x, \alpha_2^*(x) = 0. \quad (2.2.25)$$

Thus from (2.2.21) and (2.2.25) we obtain

$$\dot{x} = x, x(0) = y. \quad (2.2.26)$$

Therefore

$$x(t) = x(0)\exp(t) = y\exp(t), \quad (2.2.27)$$

and

$$J_1(\alpha_1^*) = \frac{1}{2}\int_0^\infty e^{-t}(\alpha_1^*(t))^2 dt =$$

$$= \frac{1}{2}\int_0^\infty e^{-t}x^2(t)dt = \frac{x(0)}{2}\int_0^\infty e^{-t}\exp(2t)dt = \frac{x(0)}{2}\int_0^\infty \exp(t)dt = \infty. \quad (2.2.28)$$

Thus the corresponding feedback (2.2.25) do not yield a solution to the differential game (2.2.21)-(2.2.22).

We now assume that the player have conflicting interest. Namely, their running costs $h_i(x)$ satisfy:

$$h_1'(x) \leq 0 \leq h_2'(x). \quad (2.2.29)$$

Assume that $h_i(x) = k_i x$ with $k_1 + k_2 > 0$, which is not rectrictive. The existence of an admissible solution for (2.2.18) is trivial, since we have the constant solution $p = (k_1, k_2)$, which corresponds to

$$u_1(x) = k_1 x + k_1 k_2 + \frac{1}{2}k_1^2,$$

$$u_2(x) = k_2 x + k_1 k_2 + \frac{1}{2}k_2^2. \quad (2.2.30)$$

# II.3. Infinite Horizon 2-Persons noncooperative differential game with nonlinear dinamics imbeded into a small white noise. Infinitesimal stochastic game.

Let us consider 2-persons differential game, with nonlinear dynamics

$$\dot{x} = f_1(x,\alpha_1) + f_2(x,\alpha_2(x)) = \widetilde{f}(x,\alpha_1,\alpha_2),$$

$$x(0) = y \in \mathbb{R}^n,$$

$$f_1(x,\alpha_1) = \left(f_1^1(x,\alpha_1), f_1^2(x,\alpha_1), \ldots, f_1^n(x,\alpha_1)\right), \quad (2.3.1)$$

$$f_2(x,\alpha_2) = \left(f_2^1(x,\alpha_1), f_2^2(x,\alpha_1), \ldots, f_2^n(x,\alpha_1)\right),$$

$$\alpha_j(t) \in A_j, j = 1,2$$

**Definition** 2.3.1. Let be a scalar function $V : \mathbb{R}^n \to \mathbb{R}$. $V$ is a Lyapunov-candidate-function if it is a locally positive-definite function, i.e. (i) $V(0) = 0$, (ii) $V(x) > 0, \forall x, x \in U \backslash \{0\}$, with $U$ being a neighborhood region around $x = 0$.

**Definition** 2.3.2. Let $\dot{x}_i = f_i(x_1, \ldots, x_n; \alpha_1, \alpha_2), i = 1, \ldots, n.$

$$\dot{V}(x;f) \triangleq \sum_{i=1}^{n} \frac{\partial V(x)}{\partial x_i} f(x). \quad (2.3.2)$$

**Definition** 2.3.3. Differential game (2.3.1) is dissipative iff exist Lyapunov-candidate-function $V(x)$ and constants $C > 0, R \geq 0$ such that:

$$\dot{V}(x;\widetilde{f}) \leq -CV(x), \|x\| \geq R,$$

$$V_r(x) = \lim_{r \to \infty} \left( \inf_{\|x\|>r} V(x) \right) = \infty. \tag{2.3.3}$$

$$\widetilde{f}(x,\alpha_1,\alpha_2) = \sum_{j=1}^{2} f_j(x,\alpha_j).$$

**Definition** 2.3.4. *Let us consider an 2-persons stochastic differential game with nonlinear dynamics:*

$$\frac{dx(t)}{dt} = \sum_{j=1}^{2} f_j(x,\alpha_j(x)), + \sqrt{\varepsilon}\,\dot{W}(t),$$

$$\alpha_i(t) \in A_i.$$

$$x(0) = x_0 \in \mathbb{R}^n, \tag{2.3.4}$$

$$\alpha_j(t) \in U_j, i = 1,2.$$

$$t \in [0,T], \varepsilon \ll 1.$$

**Definition** *Here $t \to \alpha_i(t)$ is the determined control chosen by the $i$-th player, within a set of admissible control values $U_i$, and the payoff for the $i$-th player:*

$$J_i(\alpha_i) \triangleq \mathbf{E}\left[\int_0^\infty e^{-t}\psi_i(x(t),\alpha_i(t))\,dt\right], \tag{2.3.5}$$

where $t \mapsto x(t,\omega)$ is the trajectory of (2.3.4).

**Definition** 2.3.5.(***Infinitesimal stochastic differential game***) *Stochastic differential game (2.3.4)-(2.3.5) is the determined differential game (2.3.1) imbeded into a 'small' white noise or infinitesimal stochastic differential game.*

# III. Infinitesimal Reduction. Strong large deviations principle for infinitesimal stochastic differential game. "Step by step" strategy.

## III.1. Infinitesimal Reduction. Strong large deviations principle for infinitesimal stochastic differential game.

Let us consider now a family $X_t^\varepsilon$ of the solutions SDE:

$$d\mathbf{X}_t^\varepsilon = \mathbf{b}(\mathbf{X}_t^\varepsilon, t) + + \sqrt{\varepsilon}\, d\mathbf{W}(t),$$

$$\mathbf{X}_t^\varepsilon(0) = x_0 \in \mathbb{R}^n, t \in [0, T], \quad (3.1.1)$$

where $\mathbf{W}(t)$ is $n$-dimensional Brownian motion and $\mathbf{b}(\circ, t) : \mathbb{R}^n \to \mathbb{R}^n$ is a polinomial transform, i.e.

$$b_i(x, t) = \sum_{|\alpha|} b_\alpha^i(x, t) x^\alpha,$$

$$\alpha = (i_1, \ldots, i_k), |\alpha| = \sum_{j=0}^k i_j, \quad (3.1.2)$$

$$i = 1, \ldots, n$$

**Definition** 3.1.1. *SDE (3.1.1) is dissipative if exist Lyapunov-candidate-function $V(x)$ and constants $C > 0, R \geq 0$ such that:*

$$\dot{V}(x;\mathbf{b}) \leq -CV(x), \|x\| \geq R,$$

$$V_r(x) = \lim_{r \to \infty} \left( \inf_{\|x\|>r} V(x) \right) = \infty. \tag{3.1.3}$$

Let us consider now a family $X_t^\varepsilon$ of the solutions dissipative SDE (3.1.1).

**Theorem** 3.1.1.( **Strong large deviations principle**).[20] For the all solutions
$\mathbf{X}_t^\varepsilon = (X_{1,t}^\varepsilon, \ldots, X_{n,t}^\varepsilon)$ dissipative SDE (3.1.1) and $\mathbb{R}$ valued parameters $\lambda_1, \ldots, \lambda_n, \lambda = (\lambda_1, \ldots, \lambda_n) \in \mathbb{R}^n$, there exists a constant $C \geq 0, C = C(T, R, \varepsilon)$ such that:

$$\liminf_{\varepsilon \to 0} \mathbf{E}\left[ \|\mathbf{X}_t^\varepsilon - \lambda\|^2 \right] \leq C \|\mathbf{U}(t, \lambda)\|^2$$

$$\lambda = (\lambda_1, \ldots, \lambda_n) \in \mathbb{R}^n \tag{3.1.4}$$

where $U(t, \lambda) = (U_1(t, \lambda), \ldots, U_n(t, \lambda))$ the solution of the linear differential master equation:

$$\frac{d\mathbf{U}(t, \lambda)}{dt} = \mathbf{J}[\mathbf{b}(\lambda, t)]\mathbf{U} + \mathbf{b}(\lambda, t),$$

$$\mathbf{U}(0, \lambda) = x_0 - \lambda, \tag{3.1.5}$$

where $\mathbf{J}[\mathbf{b}(\lambda, t)]$ the Jacobian, i.e. $\mathbf{J}$ is a $n \times n$-matrix:

$$\mathbf{J}[\mathbf{b}(\lambda,t)] = \mathbf{J}[\mathbf{b}(x,t)]|_{x=\lambda} =$$

$$= \begin{bmatrix} \dfrac{\partial b_1(x,t)}{\partial x_1} & \cdots & \dfrac{\partial b_1(x,t)}{\partial x_n} \\ \cdot & \cdots & \cdot \\ \cdot & \cdots & \cdot \\ \cdot & \cdots & \cdot \\ \dfrac{\partial b_n(x,t)}{\partial x_1} & \cdots & \dfrac{\partial b_n(x,t)}{\partial x_n} \end{bmatrix}_{x=\lambda} \quad (3.1.6)$$

**Corollary** *3.1.1. Assume the conditions of the Theorem 3.1 for any*

$$\lambda = (\lambda_1,\ldots,\lambda_n) \in \mathbb{R}^n, t \in [0,T]:$$
$$\|\mathbf{U}(t,\lambda)\| = 0 \Rightarrow \liminf_{\varepsilon \to 0} \mathbf{E}\!\left[\|\mathbf{X}_t^\varepsilon - \lambda\|^2\right] = 0 \quad (3.1.7)$$

*More precisely, for any $t \in [0,T]$ and*
$$\lambda = \lambda(t) = (\lambda_1(t),\ldots,\lambda_n(t)) \in \mathbb{R}^n$$
*sutch that*
$$U_1(t,\lambda_1(t),\ldots,\lambda_n(t)) = 0,$$

$$\cdots\cdots\cdots\cdots\cdots \quad (3.1.8)$$

$$U_n(t,\lambda_1(t),\ldots,\lambda_n(t)) = 0,$$

*the equalities is satisfaed*

$$\liminf_{\varepsilon \to 0} \mathbf{E}\!\left[\|X_{1,t}^\varepsilon - \lambda_1(t)\|^2\right] = 0,$$

$$\cdots\cdots\cdots\cdots\cdots \quad (3.1.9)$$

$$\liminf_{\varepsilon \to 0} \mathbf{E}\!\left[\|X_{n,t}^\varepsilon - \lambda_n(t)\|^2\right] = 0.$$

*Let us consider an $m$-persons stochastic differential game $SDG_{m;\bar{t}}^\delta(f,0,y_1)$, with nonlinear dynamics:*

$$\frac{dx_i(t)}{dt} = f_i(x_1,\ldots,x_m;\alpha_1,\ldots,\alpha_m) + \sqrt{\varepsilon}\,\dot{W}(t),$$

$$x(\bar{t}) = x_{\bar{t}} \in \mathbb{R}^n,$$

$$\alpha_i(t) \in U_i, i = 1,\ldots,m. \qquad (3.1.10)$$

$$t \in [\bar{t},\bar{t}+\delta],$$

$$\delta,\varepsilon \ll 1.$$

Here

$$t \to \alpha(t) = (\alpha_1[x_1(t),\ldots,x_m(t)],\ldots,\alpha_n[x_1(t),\ldots,x_m(t)]) \qquad (3.1.11)$$

*is the determined feedback control chosen by the i-th player, within a set of admissible control values $U_i$, and the payoff for the i-th player:*

$$\bar{J}_i = \mathbf{E}\left[\sum_{i=1}^{m}[x_i(\bar{t}+\delta,\omega) - y_i]^2\right] \qquad (3.1.12)$$

where $t \mapsto x(t,\omega)$ is the trajectory of the Eq. (3.1.10).

**Definition** 3.1.2. *Stochastic differential game $SDG_{m;\bar{t}}^{\delta}(f,0,y_1,y_2)$, is the infinitesimal stochastic differential game.*

Suppose that:
$$\check{\alpha}(t,\lambda;x_1(t),\ldots,x_m(t)) = \mathbf{Q}(t) \times [x_1(t),\ldots,x_m(t)]^{\mathsf{T}}, t \in [\bar{t},\bar{t}+\delta]$$

where $\mathbf{Q}(t)$ is a $m \times m$ matrix. Thus:

$$\frac{dx_i(t)}{dt} = f_i(x_1,\ldots,x_m;\breve{\alpha}_1(t,\lambda),\ldots,\breve{\alpha}_m(t,\lambda)) + \sqrt{\varepsilon}\,\dot{W}(t),$$

$$x(\bar{t}) = x_{\bar{t}} \in \mathbb{R}^n,$$

$$\breve{\alpha}_i(t) \in U_i, i = 1,\ldots,m.$$

$$t \in [\bar{t}, \bar{t}+\delta], \quad (3.1.13)$$

$$\delta, \varepsilon \ll 1.$$

$$\bar{\mathbf{J}}_i = \mathbf{E}\left[\sum_{i=1}^m [x_i(\bar{t}+\delta,\omega) - \lambda_i]^2\right].$$

*Strong large deviations principle for infinitesimal stochastic differential game.*

**Theorem.** 3.1.2. For the all solutions

$$\{\mathbf{X}_t^\varepsilon, \breve{\alpha}(t,\lambda)\} = (X_{1,t}^\varepsilon,\ldots,X_{m,t}^\varepsilon), (\breve{\alpha}_1(t,\lambda),\ldots,\breve{\alpha}_m(t,\lambda))$$

dissipative infinitesimal SDG (3.13) with $\mathbb{R}$ valued parameters $\lambda = (\lambda_1,\ldots,\lambda_m) \in \mathbb{R}^m$, there exists a constant $C \geq 0, C = C(T,R,\varepsilon)$ such that:

$$\liminf_{\varepsilon \to 0} \mathbf{E}\left[\|\mathbf{X}_t^\varepsilon - \lambda\|^2\right] \leq C\|\mathbf{U}(t,\lambda)\|^2,$$

$$\lambda = (\lambda_1,\ldots,\lambda_n) \in \mathbb{R}^n, \quad (3.1.14)$$

$$t \in [\bar{t}, \bar{t}+\delta],$$

where $\mathbf{U}(t,\lambda) = (U_1(t,\lambda),\ldots,U_n(t,\lambda))$ is the trajectory of the differential *master game:*

$$\frac{d\mathbf{U}(t,\lambda)}{dt} = \mathbf{J}[\mathbf{f}(\lambda,\check{\alpha}(t,\lambda))]\mathbf{U} + \mathbf{f}(\lambda,\check{\alpha}(t,\lambda)),$$

$$\mathbf{U}(\bar{t},\lambda) = x_{\bar{t}},$$

$$\mathbf{J}_i = \|\mathbf{U}(\bar{t}+\delta)\|^2 = \sum_{i=1}^{m}[U_i(\bar{t}+\delta)]^2.$$

(3.1.15)

Where $\mathbf{J}[\mathbf{f}(\lambda,\check{\alpha}(t,\lambda))]$ the Jacobian, i.e. $\mathbf{J}$ is a $n \times n$ -matrix:

$$\mathbf{J}[\mathbf{f}(\lambda,\check{\alpha}(t,\lambda))] = \mathbf{J}[\mathbf{f}(x,\check{\alpha}(t,x))]|_{x=\lambda} =$$

$$= \begin{bmatrix} \frac{\partial f_1(x,\check{\alpha}(t,x))}{\partial x_1} & \cdots & \frac{\partial f_1(x,\check{\alpha}(t,x))}{\partial x_n} \\ \cdot & \cdots & \cdot \\ \cdot & \cdots & \cdot \\ \cdot & \cdots & \cdot \\ \frac{\partial f_n(x,\check{\alpha}(t,x))}{\partial x_1} & \cdots & \frac{\partial f_n(x,\check{\alpha}(t,x))}{\partial x_n} \end{bmatrix}\Bigg|_{x=\lambda}$$

(3.1.16)

**Corollary** *3.1.2. Assume the conditions of the Theorem 3.1.2. Then for any $\lambda = (\lambda_1,\ldots,\lambda_n) \in \mathbb{R}^n, t \in [\bar{t},\bar{t}+\delta]$ :*

$$\|\mathbf{U}(t,\lambda)\| = 0 \Rightarrow$$

$$\liminf_{\varepsilon \to 0} \mathbf{E}\left[\|\mathbf{X}_t^\varepsilon - \lambda\|^2\right] = 0.$$

(3.1.17)

$$\liminf_{\varepsilon \to 0}\left[\min_{\alpha_i(t)}\left(\max_{\alpha_j(t),j\neq i}\mathbf{J}_i\right)\right] = 0.$$

*More precisely, for any $t \in [\bar{t},\bar{t}+\delta]$ and*
$$\lambda = \lambda(t) = (\lambda_1(t),\ldots,\lambda_n(t)) \in \mathbb{R}^n$$

*sutch that*

$$U_1(t, \lambda_1(t), \ldots, \lambda_n(t)) = 0,$$

$$\cdots\cdots\cdots\cdots\cdots\cdots \quad (3.1.18)$$

$$U_n(t, \lambda_1(t), \ldots, \lambda_n(t)) = 0,$$

*the equalities is satisfaed:*

$$\liminf_{\varepsilon \to 0} \mathbf{E}\left[\|X_{1,t}^\varepsilon - \lambda_1(t)\|^2\right] = 0,$$

$$\cdots\cdots\cdots\cdots\cdots\cdots$$

$$\liminf_{\varepsilon \to 0} \mathbf{E}\left[\|X_{n,t}^\varepsilon - \lambda_n(t)\|^2\right] = 0, \quad (3.1.19)$$

$$\liminf_{\varepsilon \to 0} \left[\min_{\alpha_i(t)} \left(\max_{\alpha_j(t), j \neq i} \mathbf{J}_i\right)\right] = 0.$$

# IV. The solution for infinite horizon 2-persons noncooperative differential game with dissipative nonlinear dinamics.

## IV.1. The solution for infinite horizon 2-persons 1-dimensional noncooperative differential game with

# dissipative nonlinear dinamics.Delay-dependent output-feedback control."Step by step"- strategy.

In this section we consider a scalar 2-persons differential game, with nonlinear dynamics

$$\dot{x}(t) = -\gamma x^3(t) + \kappa x^2(t) + \alpha_1(t) + \alpha_2(t),$$

$$0 < \gamma, 0 < \kappa, \tag{4.1.1}$$

$$x(0) = y \in \mathbb{R}.$$

The functions $t \mapsto \alpha_i(t), i = 1,2,$ represent the controls implemented by the $i$-th player, chosen within a compact set of admissible controls $\mathcal{A}_i \subset \mathbb{R}$. The game takes place on $[0,+\infty)$ and each player is subject to a running cost, exponentially discounted, of the following form

$$\mathbf{J}_i(x(t), \alpha_i(t)) \triangleq \int_0^\infty e^{-t}\left[ h_i(x(t)) + \frac{\alpha_i^2(t)}{2} \right] dt. \tag{4.1.2}$$

From Eqs.(3.15)-(3.16) we obtain master game with linear dynamics for the optimal control problem (4.1.1)-(4.1.2).

$$\dot{z} = (-3\gamma\lambda^2 + 2\kappa\lambda)z - \gamma\lambda^3 + \kappa\lambda^2 + \alpha_1(t) + \alpha_2(t),$$

$$z(0) = y - \lambda, \tag{4.1.3}$$

$$\mathbf{J}_i(z(t), \alpha_i(t)) \triangleq \int_0^\infty e^{-t}\left[ h_i(z(t)) + \frac{\alpha_i^2(t)}{2} \right] dt.$$

Set

$$z(t) + \lambda = w(t),$$

$$z(t) = w(t) - \lambda.$$

(4.1.4)

From (4.1.1) and (4.1.4) we deduce master game with linear dynamics for the optimal control problem

$$\dot{w} = (-3\gamma\lambda^2 + 2\kappa\lambda)w + 2\gamma\lambda^3 - \kappa\lambda^2 + \alpha_1(t) + \alpha_2(t),$$

$$w(0) = y,$$

$$\mathbf{J}_1(w(t), \beta_1(t)) \triangleq \int_0^\infty e^{-t}\left[h_1(w(t)) + \frac{\beta_1^2(t)}{2}\right]dt,$$

(4.1.5)

$$\beta_1(t) = \alpha_1(t) + 2\gamma\lambda^3 - \kappa\lambda^2,$$

$$\mathbf{J}_2(w(t), \alpha_2(t)) \triangleq \int_0^\infty e^{-t}\left[h_2(w(t)) + \frac{\alpha_2^2(t)}{2}\right]dt.$$

From Eq.(A.1.1) (see Appendix A) and (4.1.5) we deduce master game with simple linear dynamics for the optimal control problem:

$$\dot{w}_1(t; t_f) = \beta_1(t) \exp[(-3\gamma\lambda^2 + 2\kappa\lambda)(t_f - t)] +$$
$$+ \alpha_2(t) \exp[(-3\gamma\lambda^2 + 2\kappa\lambda)(t_f - t)],$$

$$\dot{w}_1(t; t_f) = \widetilde{\alpha}_1(t) + \widetilde{\alpha}_2(t),$$

$$w_1(t; t_f) = w(t) \exp[(-3\gamma\lambda^2 + 2\kappa\lambda)(t_f - t)],$$

$$|(-3\gamma\lambda^2 + 2\kappa\lambda)|(t_f - t) \ll 1, |(-3\gamma\lambda^2 + 2\kappa\lambda)t_f| \ll 1,$$

$$w_1(0; t_f) = w(0) \exp[(-3\gamma\lambda^2 + 2\kappa\lambda)t_f] = y \exp[(-3\gamma\lambda^2 + 2\kappa\lambda)t_f] \simeq y,$$

$$\mathbf{J}_1(w(t), \widetilde{\alpha}_1(t)) \simeq \sum_{n=0}^{\infty} \int_{t_n}^{t_{n+1}} e^{-t} \left[ h_1(w_1(t_{n+1}; t)) + \frac{1}{2}\widetilde{\alpha}_1^2(w_1(t_{n+1}; t)) \right] dt =$$

$$= \sum_{n=0}^{\infty} \int_{t_n}^{t_{n+1}} e^{-t} [h_1(w(t) \exp[(-3\gamma\lambda^2 + 2\kappa\lambda)(t_{n+1} - t)]) +$$

(4.1.6)

$$+ \frac{1}{2}\widetilde{\alpha}_1^2(w(t) \exp[(-3\gamma\lambda^2 + 2\kappa\lambda)(t_{n+1} - t)]) \Big] dt,$$

$$\mathbf{J}_2(w(t), \widetilde{\alpha}_2(t)) \simeq \sum_{n=0}^{\infty} \int_{t_n}^{t_{n+1}} e^{-t} \left[ h_2(w_1(t_{n+1}; t)) + \frac{1}{2}\widetilde{\alpha}_2^2(w_1(t_{n+1}; t)) \right] dt =$$

$$= \sum_{n=0}^{\infty} \int_{t_n}^{t_{n+1}} e^{-t} [h_2(w(t) \exp[(-3\gamma\lambda^2 + 2\kappa\lambda)(t_{n+1} - t)]) +$$

$$\frac{1}{2}\widetilde{\alpha}_2^2(w(t) \exp[(-3\gamma\lambda^2 + 2\kappa\lambda)(t_{n+1} - t)]) \Big] dt.$$

**Definition** *4.1.1. Cutting function $\Theta_\tau(t)$ :*

$$\theta_\tau(t) \triangleq \tau - t,$$

$$\eta_\tau(t) \triangleq t - \left(\textbf{ceil}\left(\frac{t}{\tau}\right) - 1\right) \cdot \tau$$

(4.1.7)

**ceil**($x$) is a part-whole number $x \in \mathbb{R}$,

$$\Theta_\tau(t) = \theta_\tau(\eta_\tau(t)).$$

From Eqs.(4.1.6)-(4.1.7) we obtain (quasy) optimal delay-dependent output-feedback control $\alpha_1^*(t - \tau; x(t - \tau))$ for the first player and (quasy) optimal delay-dependent output-feedback control $\alpha_2^*(t - \tau; x(t - \tau))$ for the second player in the next form:

$$\alpha_1^*(\tau, t; x(t - \tau)) = -u_1'(x(t - \tau)\exp[(-3\gamma x^2(t - \tau) + 2\kappa x(t - \tau))\Theta_\tau(t)]),$$

(4.1.8)

$$\alpha_2^*(\tau, t; x(t - \tau)) = -u_2'(x(t - \tau)\exp[(-3\gamma x^2(t - \tau) + 2\kappa x(t - \tau))\Theta_\tau(t)]),$$

where $u_1(x), u_2(x)$ the admissible solution for corresponding Eqs.(2.2.12). We then define the cost functions $u_i(y), i = 1, 2$:

$$u_i(y) =$$

$$=\lim_{\tau\to 0}\left[\int_0^\infty e^{-t}[h_i(x(\tau,t-\tau,y)\exp[(-3\gamma x^2(\tau,t-\tau,y)+2\kappa x(\tau,t-\tau,y))\Theta_\tau(t)])+\right.$$

$$\left.+\frac{1}{2}\alpha_i^{*2}(\tau,t;x(\tau,t-\tau,y))\right]dt\right]\simeq$$

(4.1.9)

$$\simeq \int_0^\infty e^{-t}[h_i(x(\tau,t-\tau,y)\exp[(-3\gamma x^2(\tau,t-\tau,y)+2\kappa x(\tau,t-\tau,y))\Theta_\tau(t)])+$$

$$+\frac{1}{2}u_i'^2(x(t-\tau,y)\exp[2(-3\gamma x^2(t-\tau,y)+2\kappa x(t-\tau,y))\Theta_\tau(t)])\Big]dt,$$

$$\tau \ll 1.$$

where $t \mapsto x(\tau,t,y)$ denotes the solution to the Cauchy problem

$$\dot{x}(t) = -\gamma x^3(t) + \kappa x^2(t) + \alpha_1^*(\tau,t;x(t-\tau)) + \alpha_2^*(\tau,t;x(t-\tau)),$$

(4.1.10)

$$x(0) = y.$$

**Definition** *4.1.1. Namely a solution $u$ to Eqs.(2.2.6) corresponding to differential game (4.1.1)-(4.1.2) is said to be an admissible solution if the following holds :*
*(A1) $u(x)$ is absolutely continuous and its derivative $u'(x)$ satisfies Eqs.(2.2.6) at a.e. point $x \in \mathbb{R}$.*

*(A2) The Cauchy problem:*

$$\dot{x}(t) = -\gamma x^3(t) + \kappa x^2(t) + G(x(t-\tau)),$$

$$\tau \ll 1,$$

(4.1.11)

$$x(0) = y \in \mathbb{R},$$

$$G(x(t)) = -u_1'(x(t)) - u_2'(x(t)).$$

has a globally defined solution, for every initial data $y \in \mathbb{R}$ such that $|x(t)| < C_0$.

**Theorem** 4.1.1. Consider the differential game (4.1.1)-(4.1.2) with the assumptions (A1)-(A2). Let $u_i(x), i = 1, 2$ be an admissible solution to the systems of H-J equations (2.2.6) corresponding to nonlinear differential game (4.1.1)-(4.1.2). Then the delay-dependent output-feedback controls

$$\alpha_1^*(\tau; t) = -u_1'(x(t-\tau)),$$

(4.1.12)

$$\alpha_2^*(\tau; t) = -u_2'(x(t-\tau)),$$

provide a Nash equilibrium solution in feedback form.

**Proof.** We have to show that the feedback $\alpha_i^*(\tau; t), i = 1, 2$ in (4.1.12) provides solution to the optimal feedback control problem for the $i$-th player:

$$\min_{\alpha_i(\cdot)} \int_0^\infty e^{-t}\left[h_i(x(t)) + \frac{1}{2}\alpha_i^2(\tau; t)\right]dt \qquad (4.1.13)$$

where the system has dynamics

$$\dot{x}(t) = -\gamma x^3(t) + \kappa x^2(t) + \alpha_i(\tau; t) + \alpha_j^*(\tau; t), i \neq j. \qquad (4.1.14)$$

Given an initial state $x(0) = y$, by the assumptions on $u_i$ it follows that the feedback strategy $\alpha_i(\tau; t) = \alpha_i(x(t - \tau))$ achieves a total cost given by $u_i(y)$. Now consider any absolutely continuous trajectory $t \mapsto x(\tau, t)$, with $x(0) = y$ to the Cauchy problem (4.1.10). Of course, this corresponds to the delay-dependent output-feedback control

$$\alpha_i(x(t - \tau)) = \dot{x}(t) + \gamma x^3(t) - \kappa x^2(t) - \alpha_j^*(x(t - \tau)), i \neq j. \quad (4.1.15)$$

implemented by the $i$-th player. We claim that the corresponding cost satisfies

$$\int_0^\infty e^{-t}\left[h_i(x(t)) + \frac{1}{2}(\dot{x}(t) + \gamma x^3(t) - \kappa x^2(t) - \alpha_j^*(x(t - \tau)))^2\right]dt \geq u_i(y). \quad (4.1.16)$$

To prove (4.1.16), we first observe that

$$\lim_{t \to \infty} e^{-t} u_i(x(t - \tau)) = 0. \quad (4.1.17)$$

Hence

$$u_i(y) = u_i(x(0)) =$$

$$\int_0^\infty \left[\frac{d}{dt}[e^{-t} u_i(x(t - \tau))]\right]dt. \quad (4.1.18)$$

The inequality (4.1.16) can now be established by checking that

$$e^{-t}\left[h_i(x(t)) + \frac{1}{2}(\dot{x}(t) + \gamma x^3(t) - \kappa x^2(t) - \alpha_j^*(x(t-\tau)))^2\right] \geq$$

(4.1.19)

$$e^{-t}u_i(x(t-\tau)) - e^{-t}u_i'(x(t-\tau)) \cdot \dot{x}(t-\tau).$$

Equivalently, letting $\alpha_i$ be as in (4.1.15)

$$u_i(x(t-\tau)) \leq [\alpha_i(x(t-\tau)) - u_j'x(t-\tau)]u_i'(x(t-\tau)) +$$

$$+\frac{1}{2}\alpha_i^2(x(t-\tau)) + h_i(x(t-\tau)), \qquad (4.1.20)$$

$$i \neq j.$$

This is clearly true because, by (2.2.7)

$$u_i = \min_{\alpha} \left\{\frac{1}{2}\alpha^2 + \alpha u_i' - u_i'u_j' + h_i(x)\right\}. \qquad (4.1.21)$$

**Definition** 4.1.2. *Namely a solution $\tilde{u}_i(x), i = 1, 2$ to Eqs.(2.2.6) corresponding to differential game (4.1.1)-(4.1.2) is said to be an admissible solution if the following holds :*
*(A°1) $\tilde{u}_i(x)$ is absolutely continuous and its derivative $\tilde{u}_i'(x)$*

*satisfies Eqs.(2.2.16) at a.e. point $x \in \mathbb{R}$.*

*(A°2) The Cauchy problem:*

$$\dot{x}(t) = -\gamma x^3(t) + \kappa x^2(t) + G(x(t-\tau)),$$

$$0 \leq \tau \ll 1,$$

$$x(0) = y \in \mathbb{R},$$

(4.1.21′)

$$G(x(t)) = -\tilde{u}_1'(x(t-\tau)) - \tilde{u}_2'(x(t-\tau)).$$

has a globally defined solution, for every initial data $y \in \mathbb{R}$ such that $|x(t)| < C_0$.

**Theorem** 4.1.2. Consider the differential game (4.1.1)-(4.1.2) with the assumptions (A°1)-(A°2). Let $\tilde{u}_i(x)$ be an admissible solutions to the systems of H-J equations (2.2.16) corresponding to

**Definition** linear master game (4.1.6). Then the controls

$$\alpha_1^*(x(t)) = \lim_{\tau \to 0} \alpha_1^*(\tau, t; x(\tau, t-\tau)) = -\tilde{u}_1'(x(t)),$$

(4.1.22)

$$\alpha_2^*(x(t)) = \lim_{\tau \to 0} \alpha_2^*(\tau, t; x(\tau, t-\tau)) = -\tilde{u}_2'(x(t)),$$

provide a Nash equilibrium solution in feedback form.

**Proof.** Given an initial state $x(0) = y$, from Eq.(4.1.8) we obtain

$$\alpha_1^*(x(t,y)) = \lim_{\tau \to 0} \alpha_1^*(\tau,t;x(t-\tau,y)) =$$

$$= \lim_{\tau \to 0} [-\tilde{u}_1'(x(t-\tau,y)\exp[(-3\gamma x^2(t-\tau,y)+2\kappa x(t-\tau,y))\Theta_\tau(t)])] =$$

$$= -\tilde{u}_1'(x(t,y)),$$

(4.1.23)

$$\alpha_1^*(x(t,y)) = \lim_{\tau \to 0} \alpha_2^*(\tau,t;x(t-\tau,y)) =$$

$$= \lim_{\tau \to 0} -\tilde{u}_2'(x(t-\tau,y)\exp[(-3\gamma x^2(t-\tau,y)+2\kappa x(t-\tau,y))\Theta_\tau(t)]) =$$

$$= -\tilde{u}_2'(x(t,y)).$$

From Eq.(4.1.9) we then obtain the cost functions $u_i(y), i = 1,2$ :

$$u_i(y) =$$

$$=\lim_{\tau \to 0} \left[ \int_0^\infty e^{-t}[h_i(x(\tau,t-\tau,y)\exp[(-3\gamma x^2(\tau,t-\tau,y)+2\kappa x(\tau,t-\tau,y))\Theta_\tau(t)]) + \right.$$

$$\left. + \frac{1}{2}\alpha_i^{*2}(\tau,t;x(\tau,t-\tau,y))\right]dt \Bigg] =$$

(4.1.24)

$$= \int_0^\infty e^{-t}\left[h_i(x(t,y)) + \frac{1}{2}\tilde{u}_i'^2(x(t,y))\right]dt.$$

## IV.2 Optimal control problem numerical simulation.

## "Step by step"- strategy.

Consider now the game for two players, with nonlinear dynamics

$$\dot{x} = -\gamma x^3 + \kappa x^2 + \alpha_1 + \alpha_2,$$

$$0 < \gamma, 0 < \kappa, \quad (4.2.1)$$

$$x(0) = y \in \mathbb{R}.$$

and cost functionals

$$J_i(\alpha_i) = \frac{1}{2} \int_0^\infty e^{-t} \alpha_i^2(t) dt, \quad (4.2.2)$$

$$i = 1, 2$$

From Eqs.(2.2.24) and (4.1.8) we obtain (quasy) optimal feedback control $\alpha_1^*(\tau, t; x(t-\tau))$ for the first player and (quasy) optimal feedback control $\alpha_2^*(\tau, t; x(t-\tau))$ for the second player in the next form:

$$\alpha_1^*(\tau, t; x(t-\tau)) = -x(t-\tau) \exp[(-3\gamma x^2(t-\tau) + 2\kappa x(t-\tau))\Theta_\tau(t)],$$

$$(4.2.3)$$

$$\alpha_2^*(\tau, t; x(t-\tau)) = 0.$$

Thus for the numerical simulation we obtain ODE:

$$\dot{x}(t) = -\gamma x^3(t) + \kappa x^2(t) + x(t-\tau) \exp[(-3\gamma x^2(t-\tau) + 2\kappa x(t-\tau))\Theta_\tau(t)]$$

$$(4.2.4)$$

$$x(0) = y \in \mathbb{R}.$$

From Eqs.(2.2.24) and (4.1.9) we then obtain the cost functions $u_i(y), i = 1, 2$ :

$$u_1(y) =$$

$$= \frac{1}{2} \lim_{\tau \to 0} \left[ \int_0^\infty e^{-t}[x^2(\tau,t-\tau,y)\exp[2(-3\gamma x^2(\tau,t-\tau,y)+2\kappa x(\tau,t-\tau,y))\Theta_\tau(t)]]dt \right]$$

(4.2.5)

$$u_2(y) = 0,$$

$$\tau \ll 1.$$

Where $t \mapsto x(\tau,t,y)$ denotes the solution to the Cauchy problem (4.2.4).

**Numerical simulation. Example** 4.2.1.

$\gamma = -1, \kappa = -1, y = -5, \tau = 0.01$

Below $t \mapsto z(t)$ denotes the solution to the Cauchy problem (4.2.1) for the case $\alpha_1 = \alpha_2 = 0$.

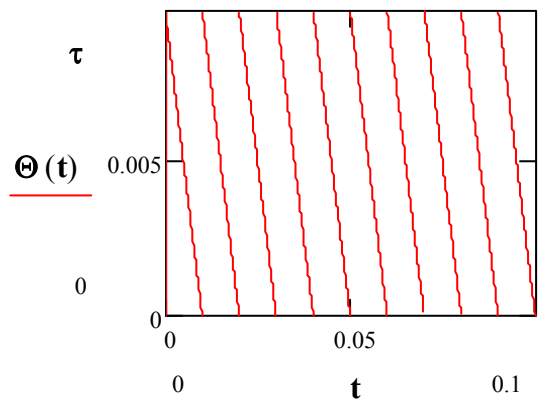

Cutting function

Pic. 4.2.1.1. $\tau = 0.01$.

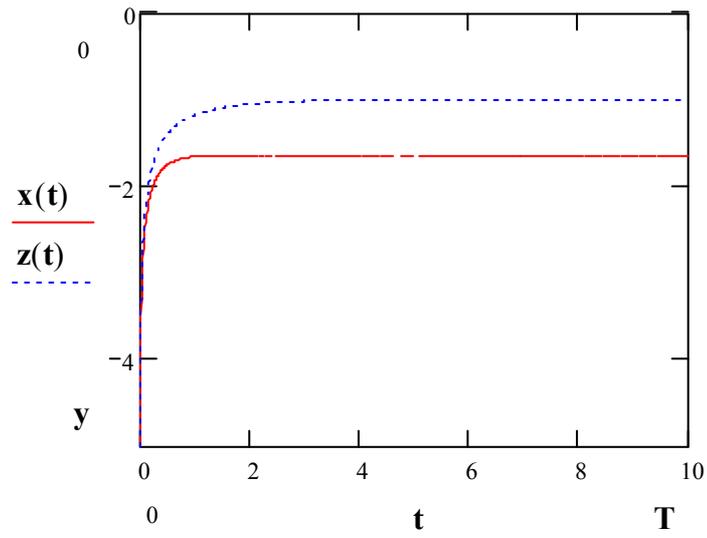

Optimal tradjectory: red curve

Pic. 4.2.1.2. $\gamma = -1, \kappa = -1, y = -5, \tau = 0.01$

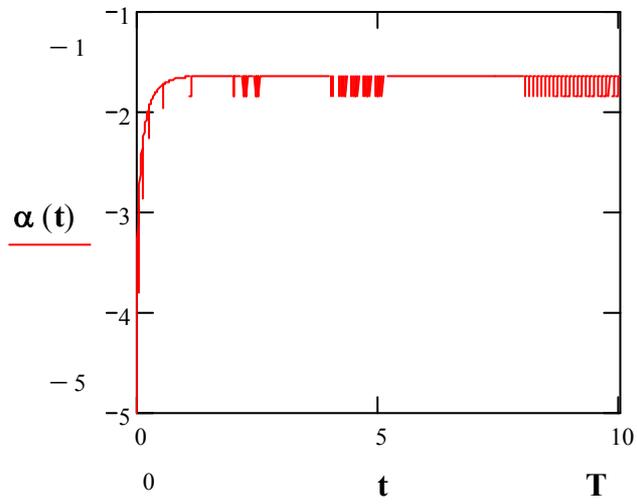

Optimal control

Pic. 4.2.1.3. $\gamma = -1, \kappa = -1, y = -5, \tau = 0.01$

$$u_1(y) = 2.64,$$
$$u_2(y) = 0.$$

Tab.4.2.1.1.Cost functions. $\gamma = -1, \kappa = -1, y = -5, \tau = 0.01$

**Numerical simulation.Example** 4.2.2.

$\gamma = -0.1, \kappa = -1, y = -5, \tau = 0.001$

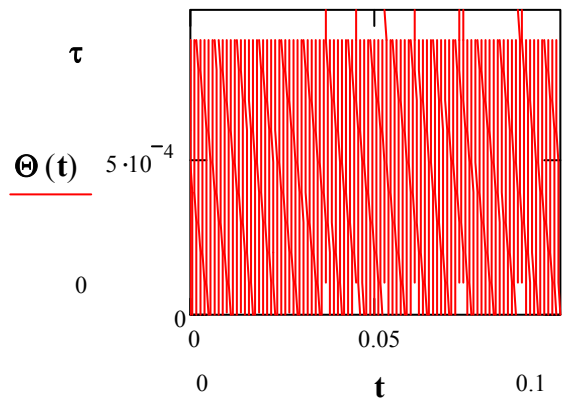

Cutting function

Pic. 4.2.2.1. $\tau = 0.001$

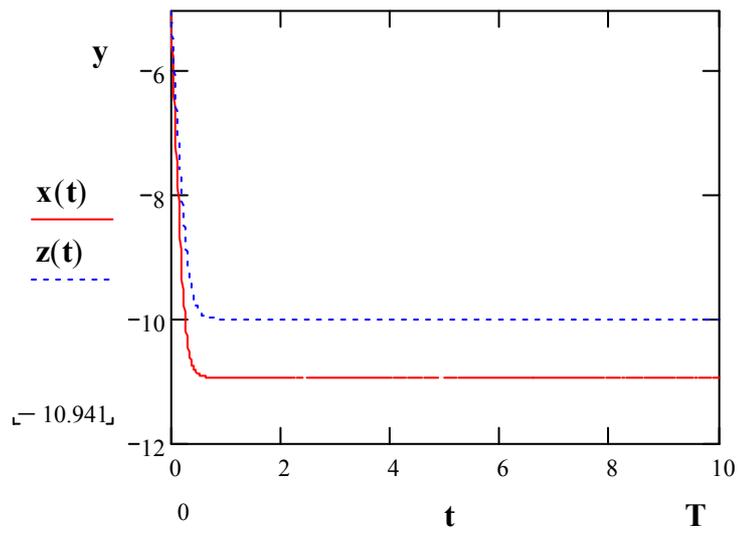

Optimal tradjectory: red curve

Pic. 4.2.2.2. $\gamma = -0.1, \kappa = -1, y = -5, \tau = 0.001$.

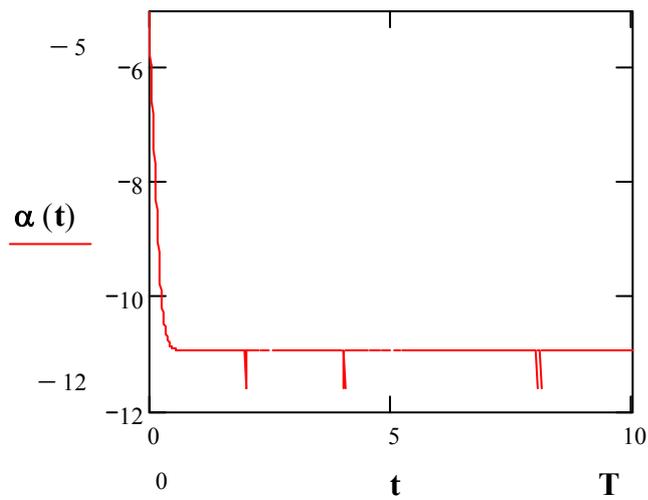

Optimal control

Pic. 4.2.2.3. $\gamma = -0.1, \kappa = -1, y = -5, \tau = 0.001$.

$$u_1(y) = 55.027,$$
$$u_2(y) = 0.$$

Tab.4.2.2.1. Cost functions. $\gamma = -0.1, \kappa = -1, y = -5, \tau = 0.001$

**Numerical simulation.Example** 4.2.3.

$\gamma = -0.1, \kappa = 1, y = -5, \tau = 0.001$

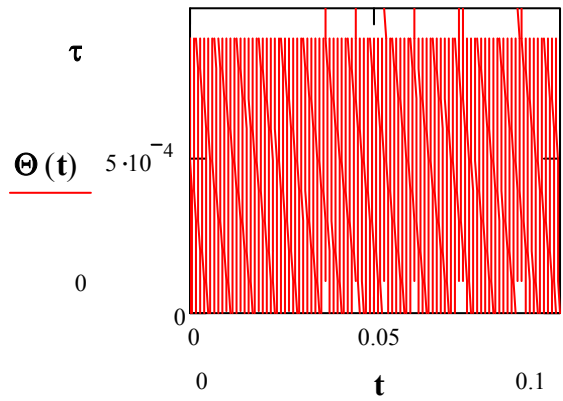

Cutting function

Pic. 4.2.3.1. $\tau = 0.001$

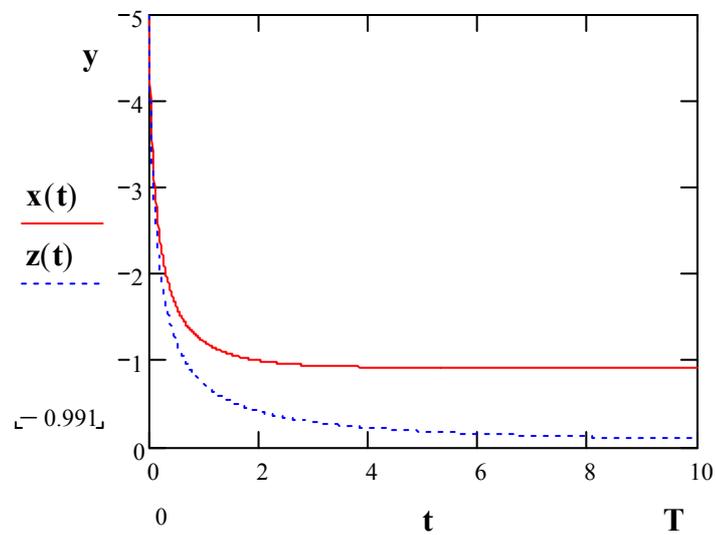

Optimal tradjectory: red curve

Pic. 4.2.3.2. $\gamma = -0.1, \kappa = 1, y = -5, \tau = 0.001$.

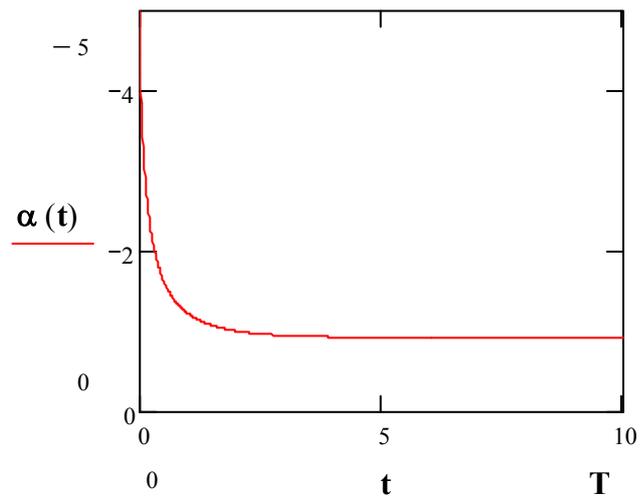

Optimal control

Pic. 4.2.3.3. $\gamma = -0.1, \kappa = 1, y = -5, \tau = 0.001$.

$$u_1(y) = 1.949,$$
$$u_2(y) = 0.$$

Tab.4.2.3.1. Cost functions. $\gamma = -0.1, \kappa = 1, y = -5, \tau = 0.001$

Consider now the game for two players, with nonlinear dynamics

$$\dot{x} = -\gamma x^3 + \kappa x^2 + \alpha_1 + \alpha_2,$$

$$0 < \gamma, 0 < \kappa, \qquad (4.2.6)$$

$$x(0) = y \in \mathbb{R}.$$

and cost functionals

$$J_i(\alpha_i) = \int_0^\infty e^{-t}\left[\frac{1}{p}(x(t)-a)^p + \frac{1}{2}\alpha_i^2(t)\right]dt,$$
(4.2.7)

$$i = 1,2$$

From Eq.(4.1.8) we obtain (quasy) optimal feedback control $\alpha_1^*(t;x(t-\tau))$ for the first player and (quasy) optimal feedback control $\alpha_2^*(t;x(t-\tau))$ for the second player in the next form:

$$\alpha_1^*(\tau,t;x(t-\tau)) = -u_1'[x(t-\tau)\exp[(-3\gamma x^2(t-\tau) + 2\kappa x(t-\tau))\Theta_\tau(t)]],$$
(4.2.8)
$$\alpha_2^*(\tau,t;x(t)) = -u_2'[x(t-\tau)\exp[(-3\gamma x^2(t-\tau) + 2\kappa x(t-\tau))\Theta_\tau(t)]].$$

Functions $u_i(x)$ satisfies Eq.(2.2.17) and its derivative $p_i(x) = u_i'(x)$ satisfies Eq.(2.2.18). Thus from (4.2.7) we obtain

$$p_1'(x) = (x-a)^{p-1}p_2(x) - p_1^2(x),$$
(4.2.9)
$$p_2'(x) = (x-a)^{p-1}p_1(x) - p_2^2(x).$$

**Numerical simulation. Example** 4.2.4.

$\gamma = -0.1, \kappa = 3, p = 3, y = 5, \tau = 0.01$

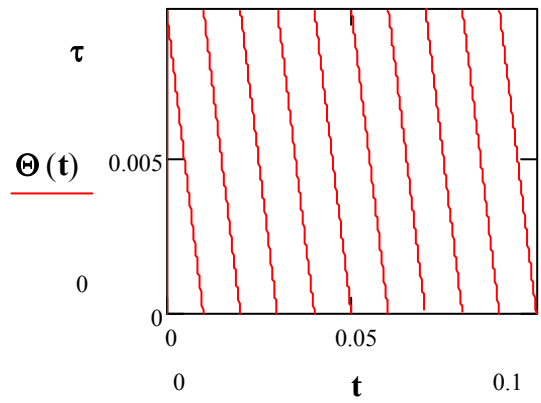

Cutting function

Pic. 4.2.4.1. $\tau = 0.01$.

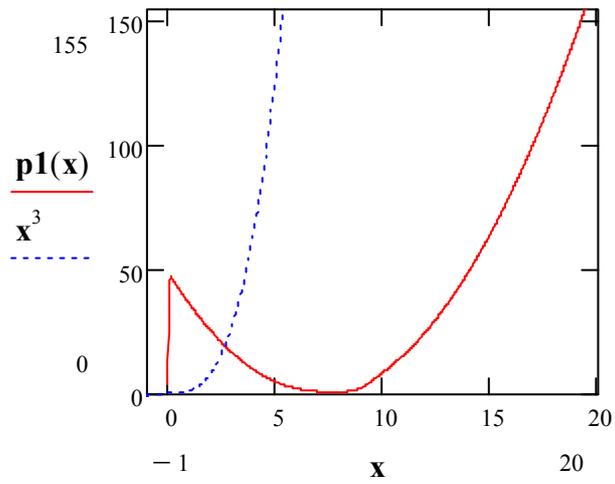

Pic. 4.2.4.2. $p_1(x)$. $p_2(0) = p_1(0) = 4, \gamma = -0.1, \kappa = 3, p = 3, y = 5, \tau = 0.01$

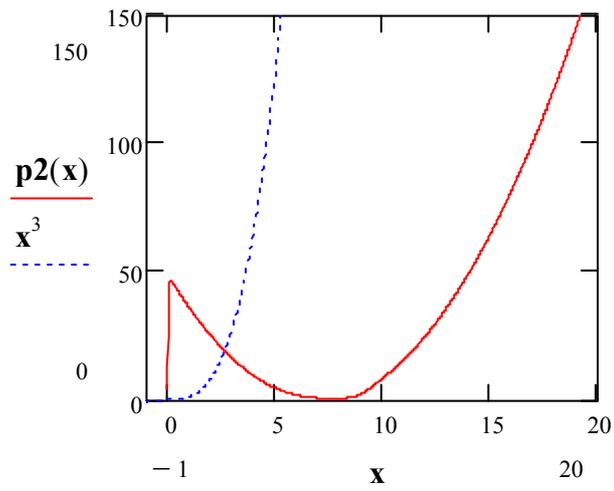

Pic. 4.2.4.3. $p_2(x)$. $p_1(0) = p_2(0) = 4, \gamma = -0.1, \kappa = 3, p = 3, y = 5, \tau = 0.01$

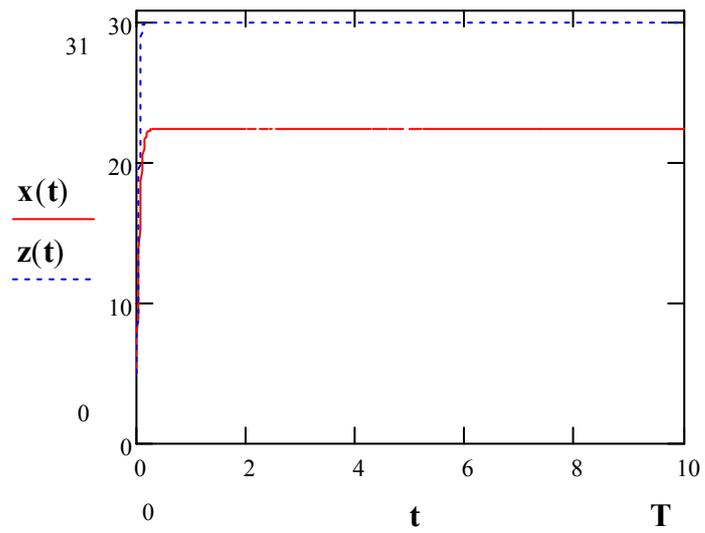

Optimal tradjectory: red curve

Pic. 4.2.4.4. $x(t)$. $p_1(0) = p_2(0) = 4, \gamma = -0.1, \kappa = 3, p = 3, y = 5, \tau = 0.01$

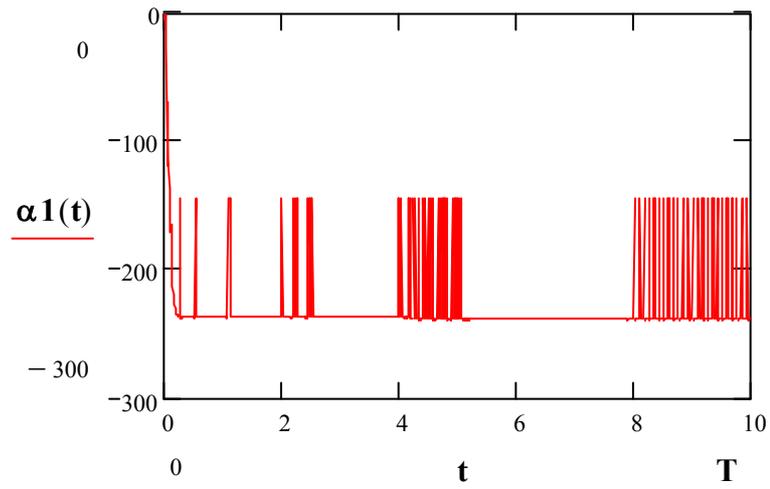

optimal feedback control for the first player

Pic. 4.2.4.5. $\alpha_1(x(t-\tau),y)$. $\gamma = -0.1, \kappa = 3, p = 3, y = 5, \tau = 0.01$

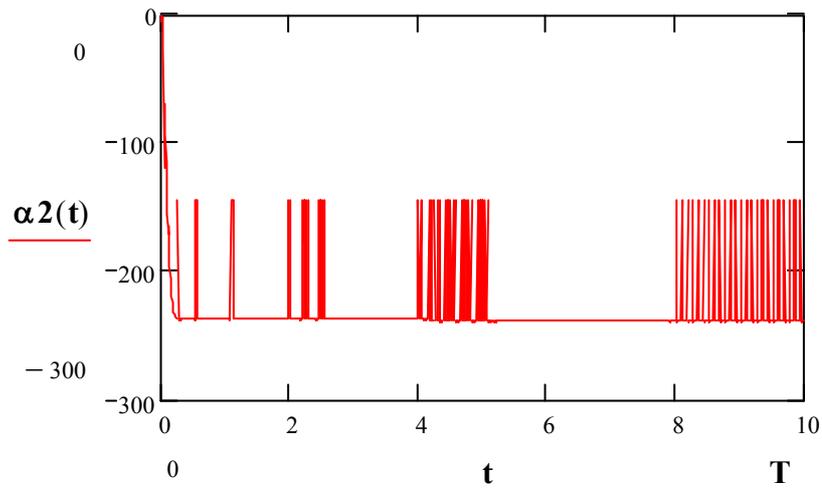

optimal feedback control for the second player

Pic. 4.2.4.6. $\alpha_2(x(t-\tau),y)$. $\gamma = -0.1, \kappa = 3, p = 3, y = 5, \tau = 0.01$

**Numerical simulation. Example** 4.2.5.

$\gamma = -0.1, \kappa = 3, p = 3, y = 5, \tau = 0.001$

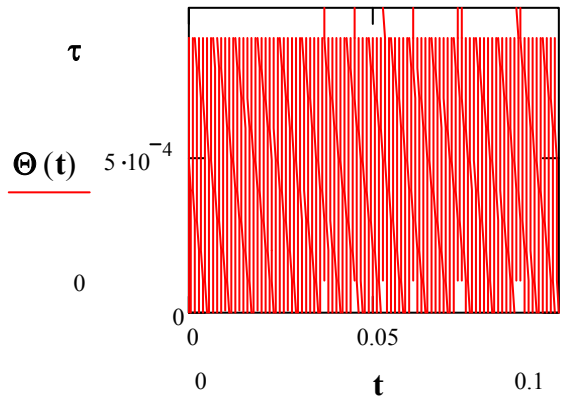

Cutting function

Pic. 4.2.5.1. $\tau = 0.001$.

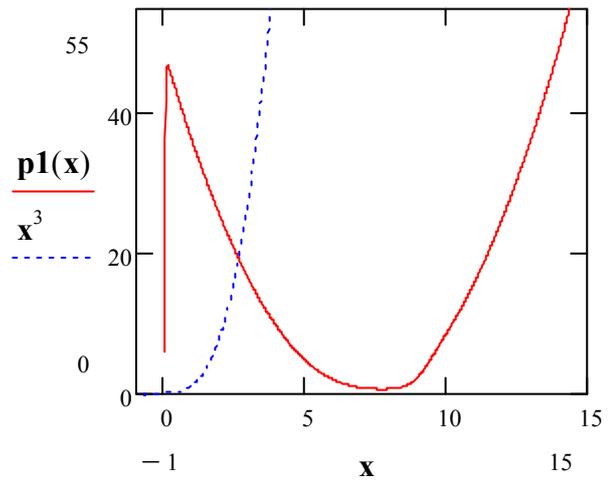

Pic. 4.2.5.2. $p_1(x)$. $p_2(0) = p_1(0) = 4, \gamma = -0.1, \kappa = 3, p = 3, y = 5, \tau = 0.001$

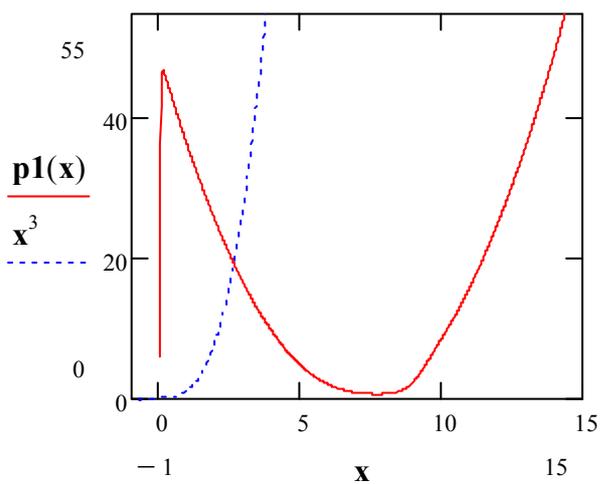

Pic. 4.2.5.3. $p_2(x)$. $p_1(0) = p_2(0) = 4, \gamma = -0.1, \kappa = 3, p = 3, y = 5, \tau = 0.001$

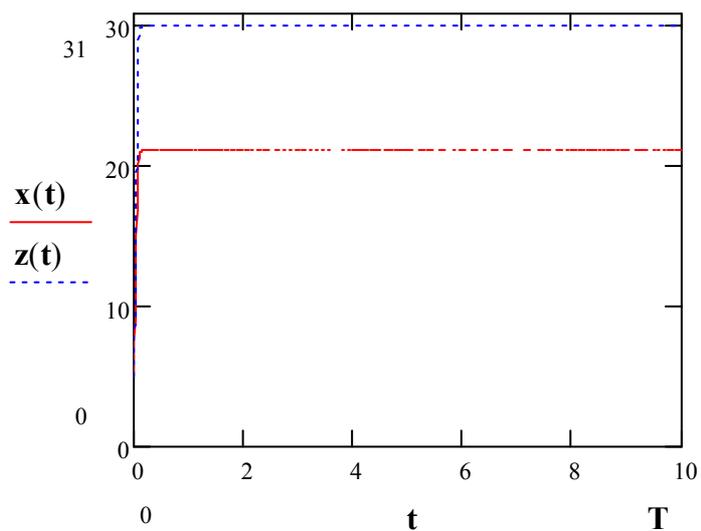

Optimal tradjectory: red curve

Pic. 4.2.5.4. $x(t)$. $p_1(0) = p_2(0) = 4, \gamma = -0.1, \kappa = 3, p = 3, y = 5, \tau = 0.001$

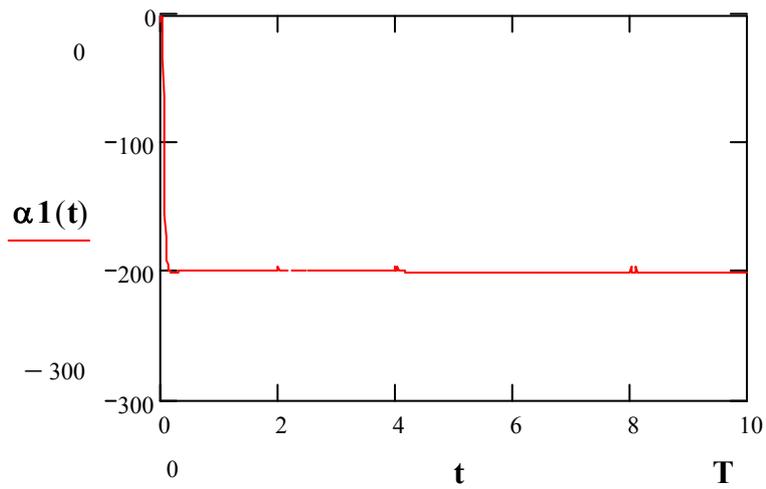

optimal feedback control for the first player

Pic. 4.2.5.5. $\alpha_1(x(t-\tau), y)$. $\gamma = -0.1, \kappa = 3, p = 3, y = 5, \tau = 0.001$

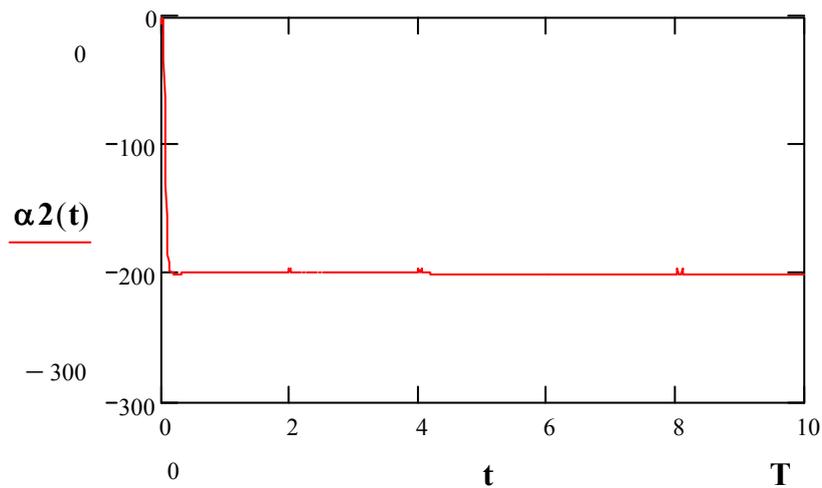

optimal feedback control for the second player

Pic. 4.2.5.6. $\alpha_2(x(t-\tau), y)$. $\gamma = -0.1, \kappa = 3, p = 3, y = 5, \tau = 0.001$

## IV.3. Players with conflicting interests.

We consider here a game for two players, with dynamics (4.1.1) and cost functionals as in (4.1.2).We assume that the player have conflicting interest. Namely,their running costs $h_i(x)$ satisfy

$$h'_1(x) \leq 0 \leq h'_2(x). \qquad (4.2.10)$$

$$\dot{x} = -\gamma x^3 + \kappa x^2 + \alpha_1 + \alpha_2,$$

$$0 < \gamma, 0 < \kappa,$$

$$x(0) = y \in \mathbb{R}.$$

$$\mathbf{J}_i(x(t), \alpha_i(t)) \triangleq \int_0^\infty e^{-t}\left[h_i(x(t)) + \frac{\alpha_i^2(t)}{2}\right]dt. \qquad (4.2.11)$$

$$h_1(x) = -\frac{1}{p}(x(t) - a_1)^p,$$

$$h_2(x) = \frac{1}{p}(x(t) - a_2)^p.$$

From Eq.(4.1.8) we obtain (quasy) optimal feedback control $\alpha_1^*(t; x(t-\tau))$

for the first player and (quasy) optimal feedback control $\alpha_2^*(t; x(t-\tau))$ for the second player in the next form:

$$\alpha_1^*(\tau, t; x(t-\tau)) = -u'_1[x(t-\tau)\exp[(-3\gamma x^2(t-\tau) + 2\kappa x(t-\tau))\Theta_\tau(t)]], \qquad (4.2.12)$$

$$\alpha_2^*(\tau, t; x(t)) = -u'_2[x(t-\tau)\exp[(-3\gamma x^2(t-\tau) + 2\kappa x(t-\tau))\Theta_\tau(t)]].$$

Functions $u_i(x)$ satisfies Eq.(2.2.17) and its derivative $p_i(x) = u'_i(x)$ satisfies Eq.(2.2.18). Thus from (4.2.7) we obtain

$$p'_1(x) = -\left[(x-a_1)^{p-1} + (x-a_2)^{p-1}\right]p_1(x) - (x-a_1)^{p-1}p_2(x) - p_1^2(x),$$

$$p'_2(x) = \left[(x-a_1)^{p-1} + (x-a_2)^{p-1}\right]p_2(x) + (x-a_2)^{p-1}p_1(x) - p_2^2(x). \qquad (4.2.13)$$

Thus for the numerical simulation we obtain ODE:

$$\dot{x} = -\gamma x^3 + \kappa x^2 + \alpha_1^*(\tau, t; x(t-\tau)) + \alpha_2^*(\tau, t; x(t-\tau))$$

(4.2.14)

$$x(0) = y \in \mathbb{R}.$$

From Eqs.(2.2.30) and (4.1.9) we then obtain the cost functions $u_i(y), i = 1, 2$ :

$$u_i(y) =$$
$$= \lim_{\tau \to 0} \left[ \int_0^\infty e^{-t} h_i(x(t,y) \exp[(-3\gamma x^2(t,y) + 2\kappa x(t,y))\Theta_\tau(t)]) + \frac{1}{2}[\alpha_i^*(\tau, t; x(t-\tau))]^2 \, dt \right]$$

(4.2.9)

$$\simeq \int_0^\infty e^{-t} \left[ h_i(x(t,y) \exp[(-3\gamma x^2(t,y) + 2\kappa x(t,y))\Theta_\tau(t)]) + \frac{1}{2}[\alpha_i^*(\tau, t; x(t-\tau))]^2 \right] dt,$$

$$\tau \ll 1.$$

Where $t \mapsto x(t - \tau, y)$ denotes the solution to the Cauchy problem (4.2.14).

**Numerical simulation.Example** 4.3.1.

$$\gamma = -0.1, \kappa = 3, p = 3, y = 5, a_1 = a_2 = 1, \tau = 0.001$$

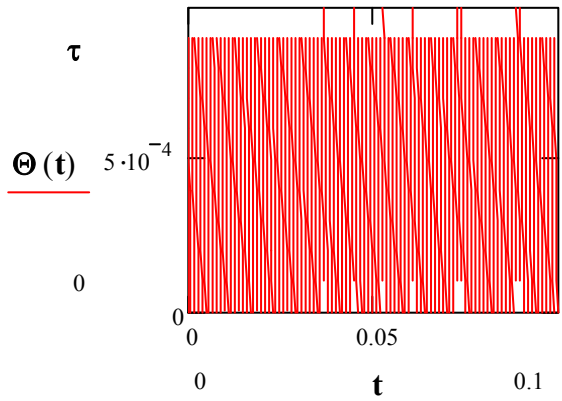

Cutting function

Pic. 4.3.1.1. $\Theta_\tau(t)$. $\tau = 0.001$.

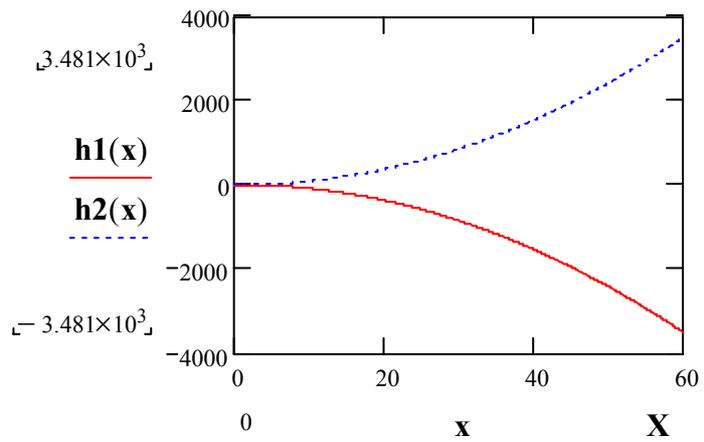

Players with conflicting interests

Pic. 4.3.1.2. $h'_1(x) = -(x - a_1)^{p-1}$, $h'_2(x) = -(x - a_2)^{p-1}$
$\gamma = -0.1, \kappa = 3, p = 3, y = 5, a_1 = a_2 = 1, \tau = 0.001$

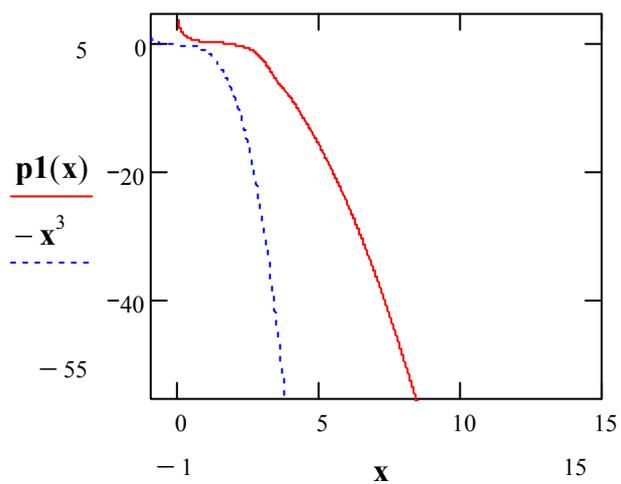

Pic. 4.3.1.3. $p_1(x)$. $p_1(0) = 4, p_2(0) = 0, \gamma = -0.1, \kappa = 3, p = 3, y = 5,$

$a_1 = a_2 = 1, \tau = 0.001.$

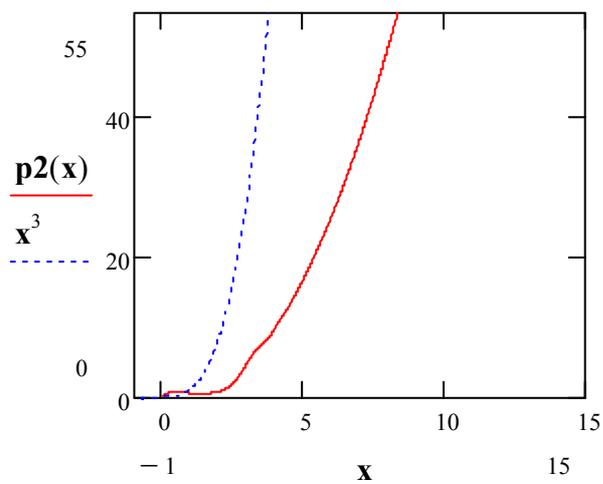

Pic. 4.3.1.4. $p_2(x)$. $p_1(0) = 4, p_2(0) = 0, \gamma = -0.1, \kappa = 3, p = 3, y = 5,$

$a_1 = a_2 = 1, \tau = 0.001.$

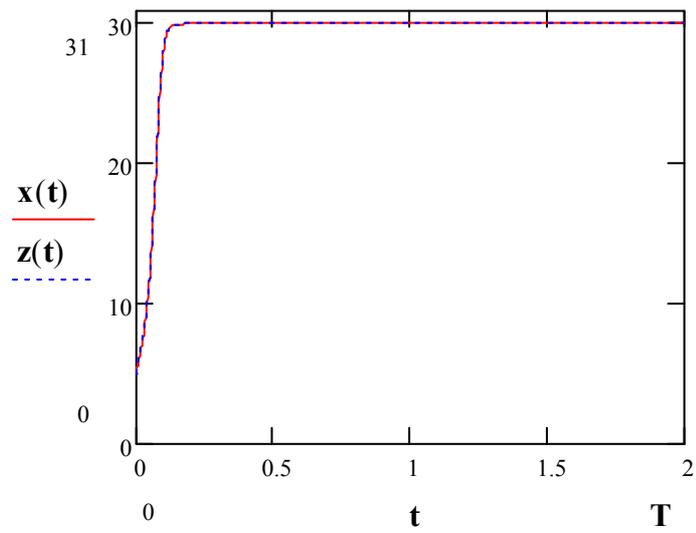

Optimal tradjectory: red curve

Pic. 4.3.5.5. $x(t)$. $\gamma = -0.1, \kappa = 3, p = 3, y = 5, p_1(0) = 4, p_2(0) = 0, a_1 = a_2 = 1, \tau = 0.001$

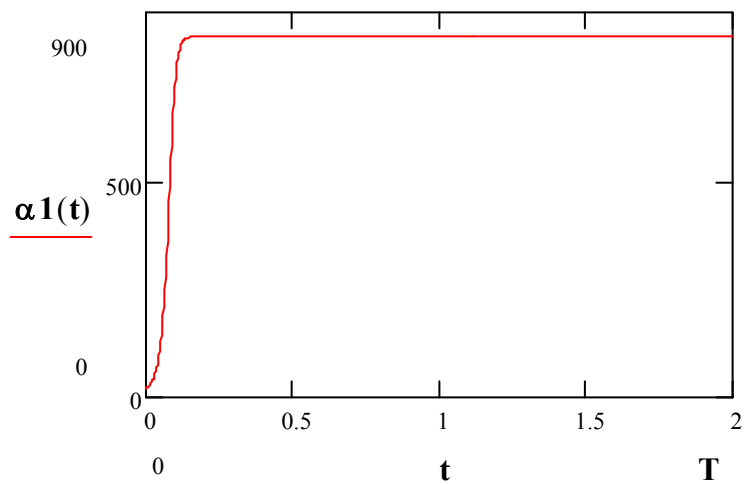

optimal feedback control for the first player

Pic. 4.3.5.6. $\alpha_2(x(t-\tau), y)$. $\gamma = -0.1, \kappa = 3, p = 3, y = 5, a_1 = a_2 = 1, \tau = 0.001$

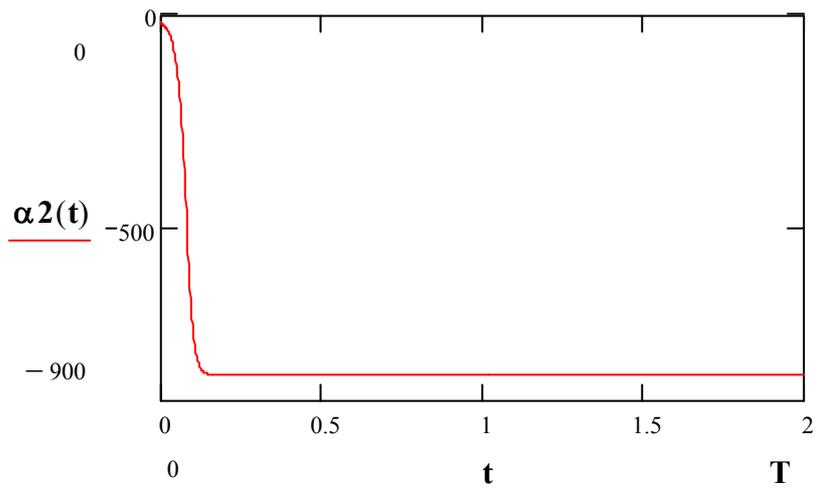

optimal feedback control for the second player

Pic. 4.3.5.7. $\alpha_2(x(t-\tau), y)$. $\gamma = -0.1, \kappa = 3, p = 3, y = 5, a_1 = a_2 = 1, \tau = 0.001$

## Numerical simulation. Example 4.3.2.

$\gamma = -0.1, \kappa = 3, p = 3, y = 5, a_1 = 1, a_2 = -3, p_1(0) = 4, p_2(0) = 0, \tau = 0.001$.

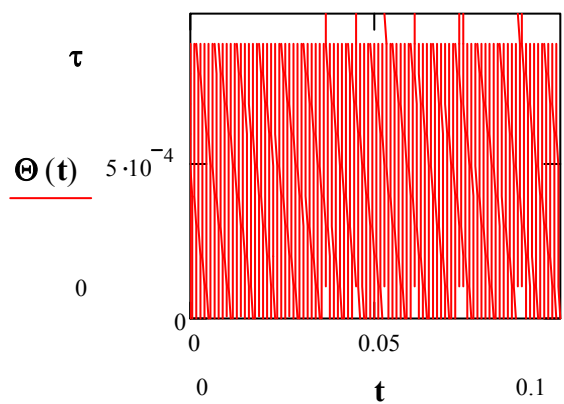

Cutting function

Pic. 4.3.2.1. $\Theta_\tau(t)$. $\tau = 0.001$

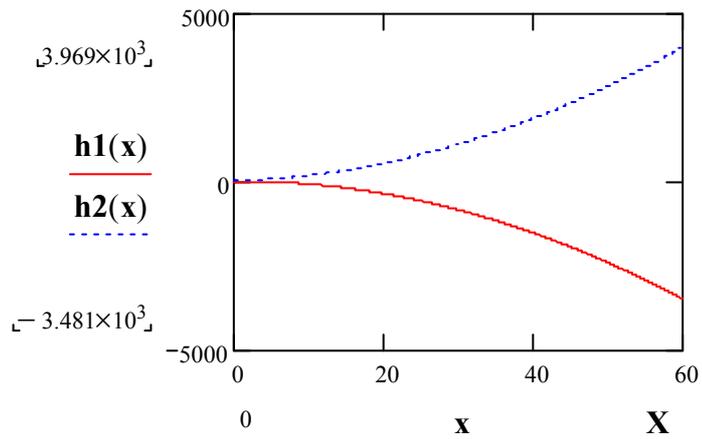

Players with conflicting interests

Pic. 4.3.2.2. $h'_1(x) = -(x - a_1)^{p-1}, h'_2(x) = -(x - a_2)^{p-1}$
$\gamma = -0.1, \kappa = 3, p = 3, y = 5, a_1 = 1, a_2 = -3, p_1(0) = 4, p_2(0) = 0, \tau = 0.001$

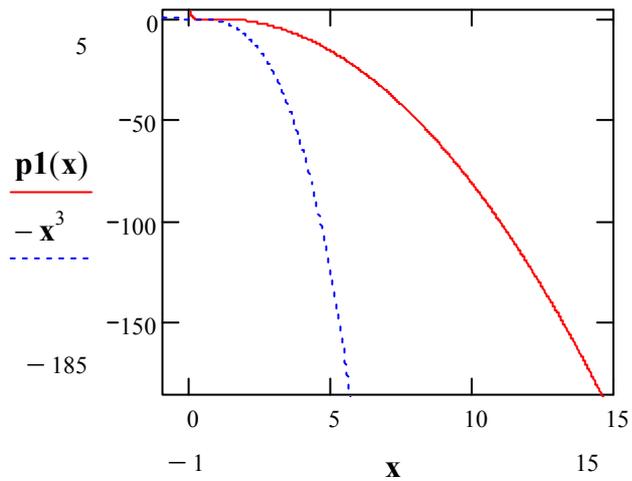

Pic. 4.3.2.3. $p_1(x). p_1(0) = 4, p_2(0) = 0, \gamma = -0.1, \kappa = 3, p = 3, y = 5,$
$a_1 = 1, a_2 = -3, \tau = 0.001.$

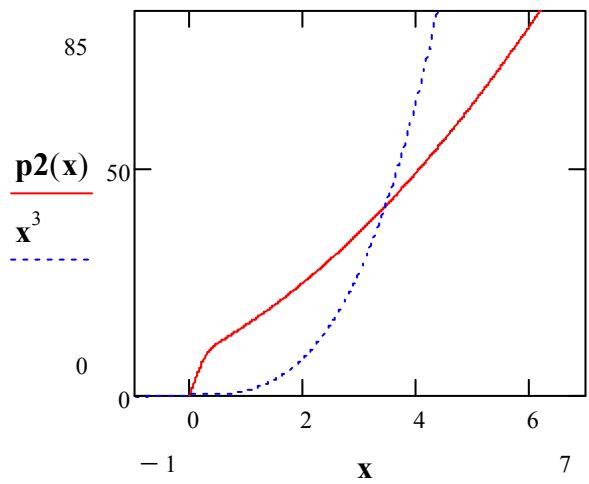

Pic. 4.3.2.4. $p_2(x)$. $p_1(0) = 4, p_2(0) = 0, \gamma = -0.1, \kappa = 3, p = 3, y = 5,$

$a_1 = 1, a_2 = -3, \tau = 0.001.$

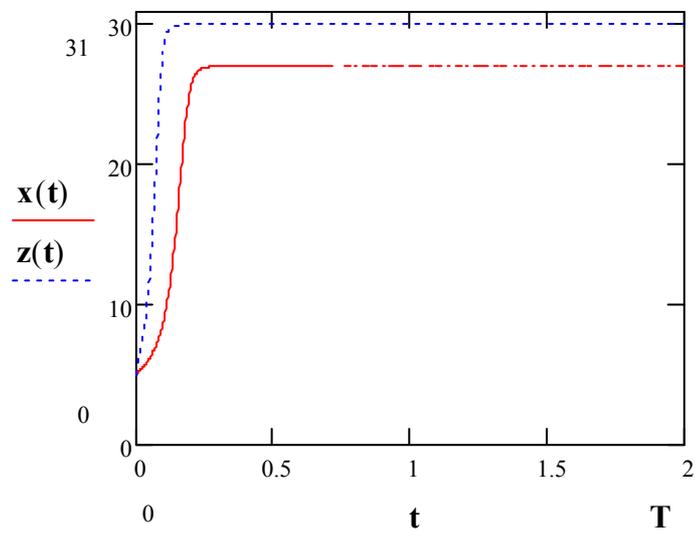

Optimal tradjectory: red curve

Pic. 4.3.2.5. $x(t)$.

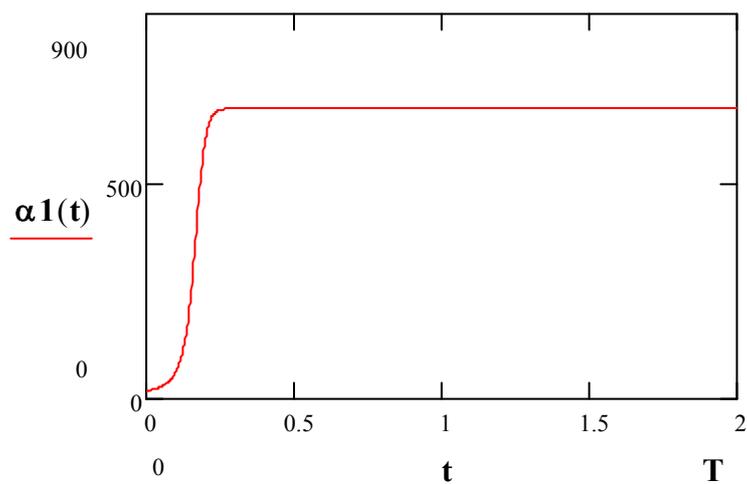

optimal feedback control of the first player

Pic. 4.3.2.6. $\alpha_1(x(t-\tau), y)$. $\gamma = -0.1, \kappa = 3, p = 3, y = 5,$

$$a_1 = 1, a_2 = -3, \tau = 0.001.$$

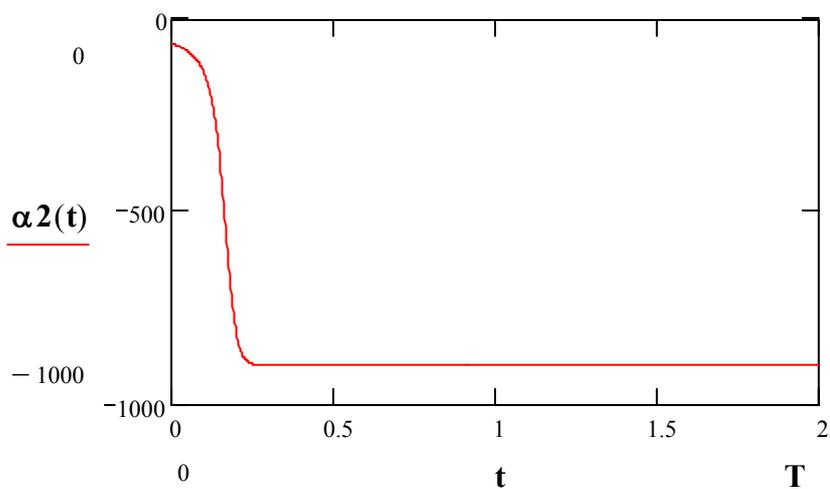

optimal feedback control of the second player

Pic. 4.3.2.7. $\alpha_2(x(t-\tau), y)$. $\gamma = -0.1, \kappa = 3, p = 3, y = 5,$

$$a_1 = 1, a_2 = -3, \tau = 0.001.$$

**Numerical simulation. Example** 4.3.3.

$\gamma = -0.05, \kappa = 3, p = 3, y = 5, a_1 = 1, a_2 = -4, p_1(0) = 4, p_2(0) = 0, \tau = 0.001.$

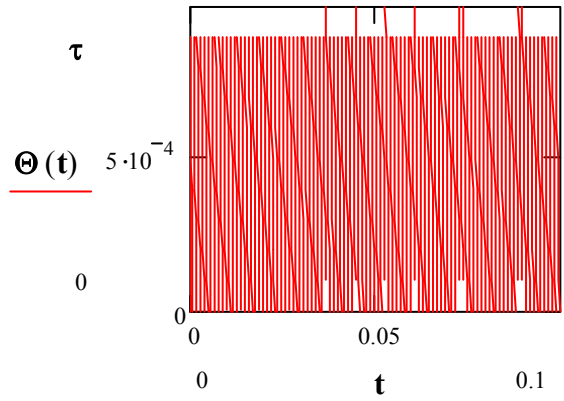

Cutting function

Pic. 4.3.3.1. $\Theta_\tau(t). \tau = 0.001$

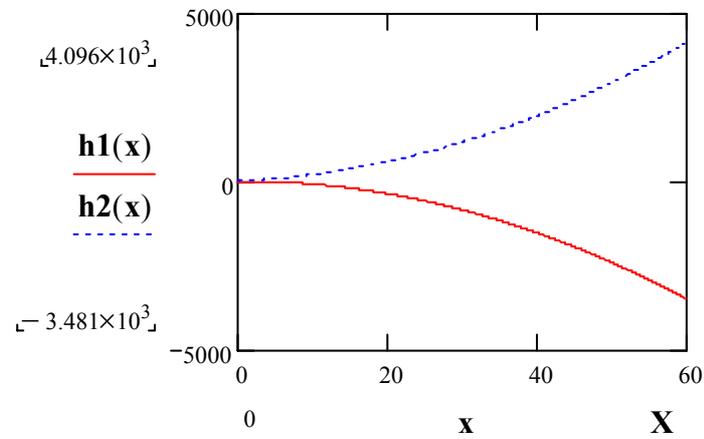

Players with conflicting interests

Pic. 4.3.3.2. $h_1'(x) = -(x - a_1)^{p-1}, h_2'(x) = -(x - a_2)^{p-1}.$

$\gamma = -0.05, \kappa = 3, p = 3, y = 5, a_1 = 1, a_2 = -4, p_1(0) = 4, p_2(0) = 0, \tau = 0.001.$

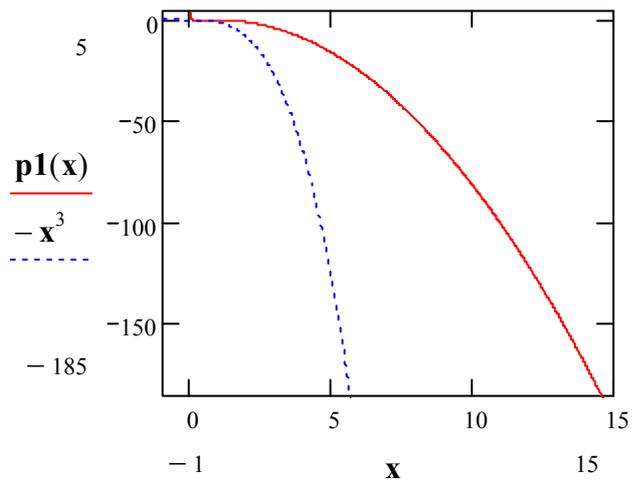

Pic. 4.3.3.3

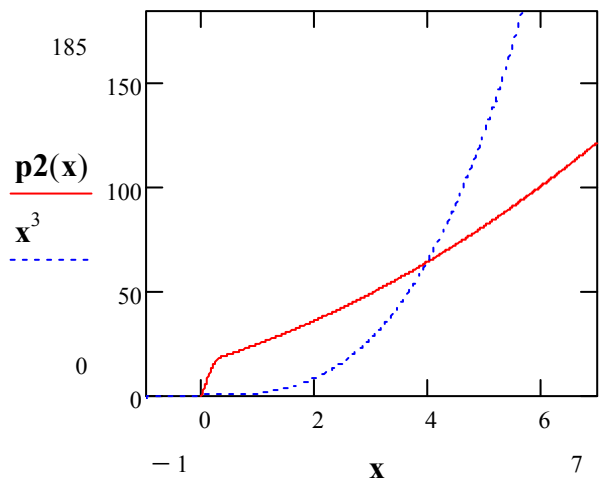

Pic. 4.3.3.4

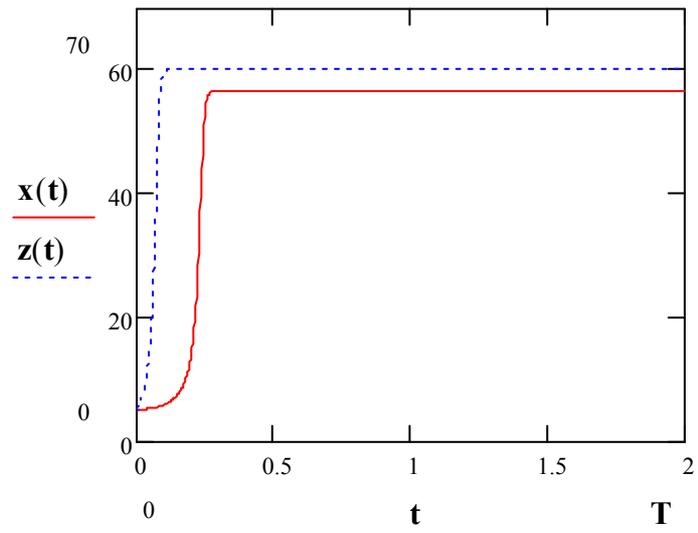

Optimal tradjectory: red curve

Pic. 4.3.3.5

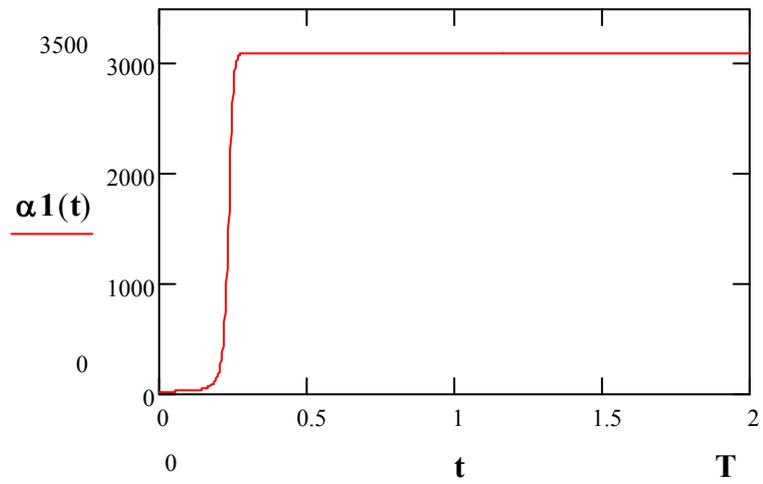

optimal feedback control of the first player

Pic. 4.3.3.6

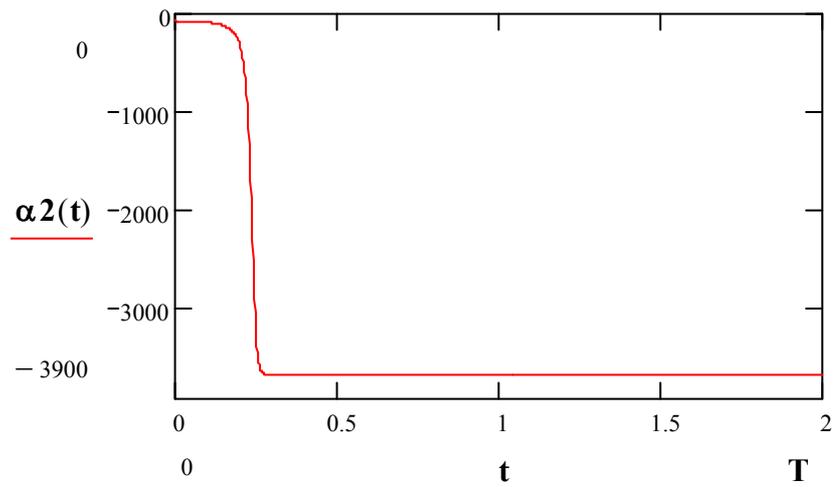

optimal feedback control of the second player

Pic. 4.3.3.7

**Numerical simulation. Example** 4.3.4.

$\gamma = -0.05, \kappa = 3, p = 3, y = 5, a_1 = 1, a_2 = -4, p_1(0) = 4, p_2(0) = 0, \tau = 0.005.$

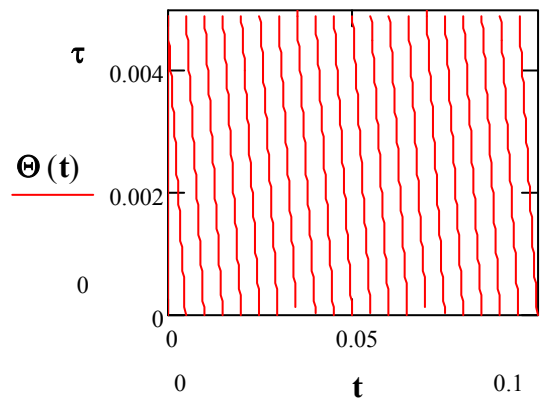

Cutting function

Pic. 4.3.4.1. $\Theta_\tau(t). \tau = 0.005$

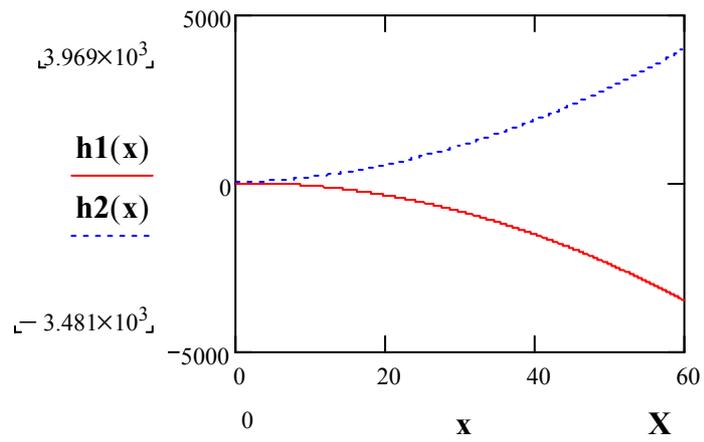

Players with conflicting interests

Pic. 4.3.4.2. $h'_1(x) = -(x-a_1)^{p-1}, h'_2(x) = -(x-a_2)^{p-1}$
$\gamma = -0.1, \kappa = 3, p = 3, y = 5, a_1 = 1, a_2 = -3, p_1(0) = 4, p_2(0) = 0, \tau = 0.001$

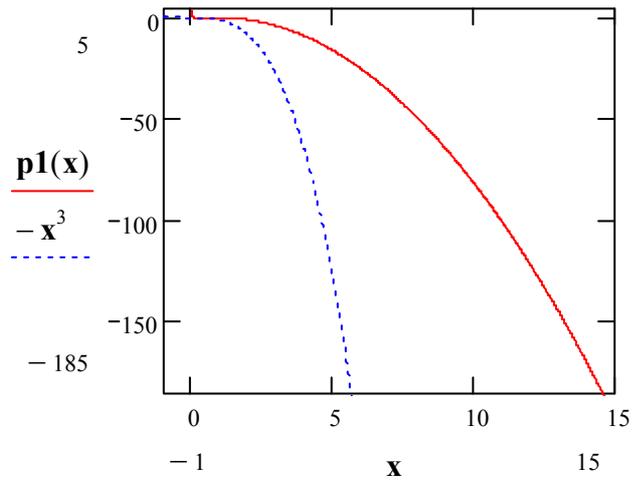

Pic. 4.3.4.3

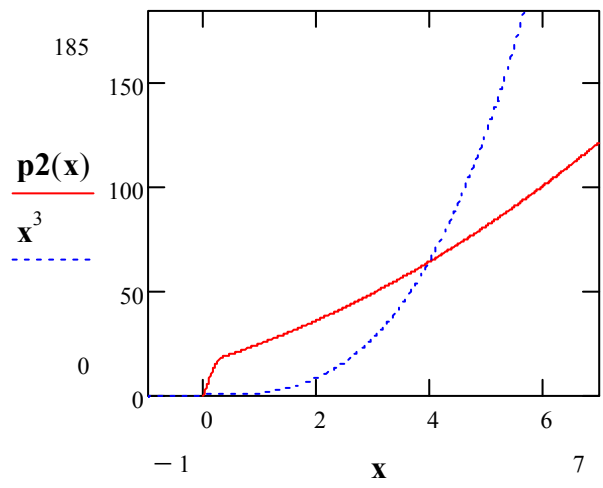

Pic. 4.3.4.4

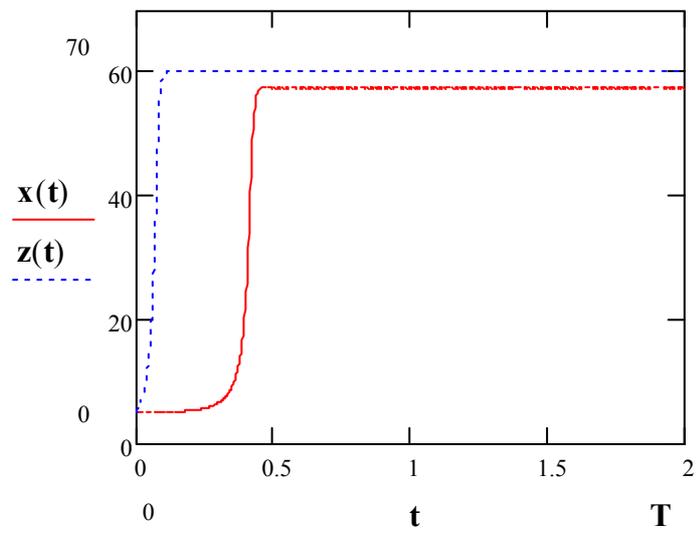

Optimal tradjectory: red curve

Pic. 4.3.4.5

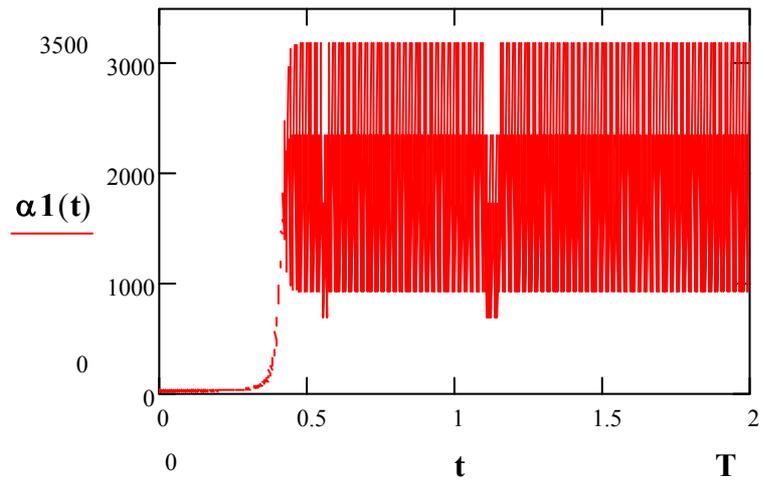

optimal feedback control of the first player

Pic. 4.3.4.6

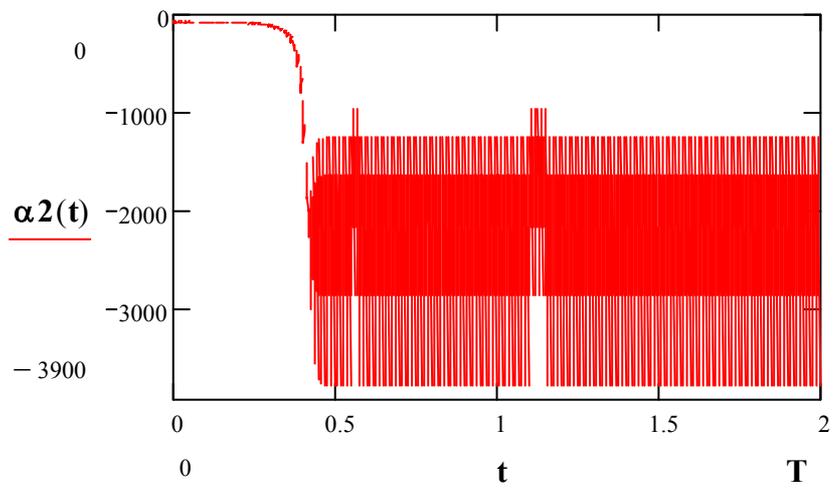

optimal feedback control of the second player

Pic. 4.3.4.7

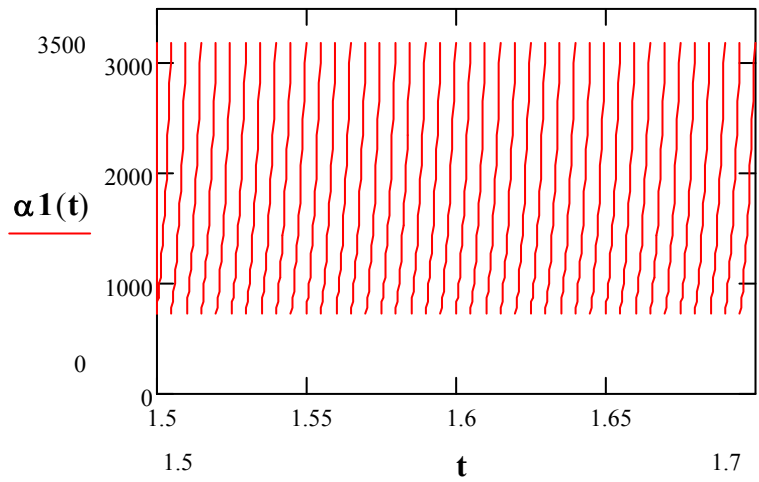

Optimal feedback control of the first player

Pic. 4.3.4.8

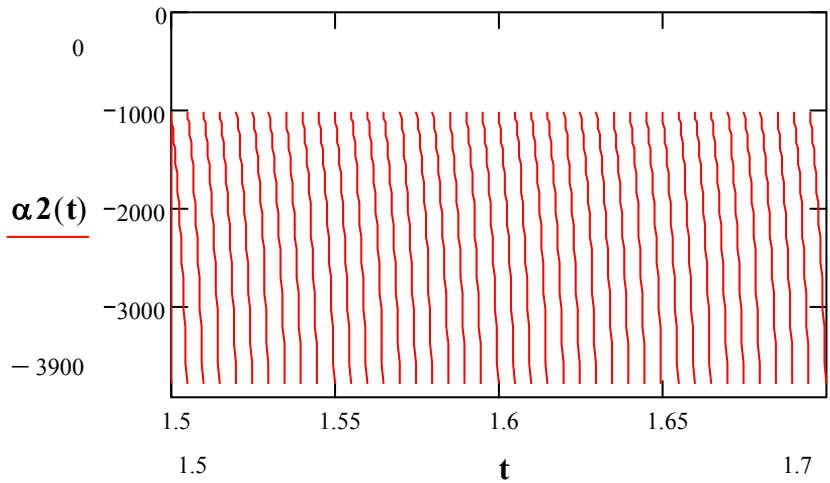

Optimal feedback control of the second player

Pic. 4.3.4.9

# Appendix.A. Reduction 2-Persons differential game

## with general linear dynamics.

In more details look a paper [20].

Let us consider an 2-persons antagonistic differential game $LDG_{2;T}(\mathbf{f}, \mathbf{0}, \mathbf{M}, \mathbf{0})$, with linear dynamics:

$$\dot{x} = A(t)x + B(t)u + C(t)v$$

$$u \in U, v \in V.  \qquad (A.1)$$

$$J = \|\mathbf{M}x(T)\|.$$

For the solutions of ordinary differential equation

$$\dot{x} = A(t)x + f(t), t \geq t_0,  \qquad (A.2)$$

the Cauchy formula

$$x(t) = X(t)x(t_0) + \int_{t_0}^{t} \Phi(t,s) \cdot f(s)ds  \qquad (A.3)$$

$X(t) \triangleq n \times n$ matrix whose columns constitute $n$ lynearly independet solutions to ordinary differential equation

$$\dot{x}(t) = A(t)x(t),  \qquad (A.4)$$

is the fundamental matrix.

$$\dot{X}(t) = A(t)X(t),$$

$$\Phi(t,s) = X(t)X^{-1}(s). \tag{A.5}$$

**Notation** A.1. $\Phi(T,t) \triangleq$ *Cauchy matrix for the solutions of ordinary differential equation (A.4).*

Thus

$$\dot{\Phi}(T,t) = \frac{\partial \Phi(T,t)}{\partial t} = -\Phi(T,t)A(t). \tag{A.6}$$

**Notation** A.2.

$$y \triangleq \boldsymbol{M}\Phi(T,t)x,$$

$$X(T,t) \triangleq \boldsymbol{M}\Phi(T,t)B(t), \tag{A.7}$$

$$Y(T,t) \triangleq \boldsymbol{M}\Phi(T,t)C(t).$$

Thus

$$\dot{y}(T,t) = \boldsymbol{M}\dot{\Phi}(T,t)x(t) + \boldsymbol{M}\Phi(T,t)\dot{x}(t). \tag{A.8}$$

Substitution (A.1) and (A.7) into (A.8) gives:

$$\dot{y}(T,t) = \mathbf{M}\dot{\Phi}(T,t)x(t) + \mathbf{M}\Phi(T,t)[A(t)x(t) + B(t)u(t) + C(t)v(t)] =$$

$$= [\mathbf{M}\dot{\Phi}(T,t) + \mathbf{M}\Phi(T,t)A(t)]x(t) + \mathbf{M}\Phi(T,t)B(t)u(t) + \mathbf{M}\Phi(T,t)C(t)v(t) =$$

$$X(T,t)u(t) + Y(T,t)v(t), i.e. \tag{A.9}$$

$$\dot{y}(T,t) = X(T,t)u(t) + Y(T,t)v(t).$$

$$J = \|y(T)\|.$$

Thus from 2-persons differential game $DG_{2;T}(\mathbf{f},\mathbf{0},\mathbf{M},\mathbf{0})$ (A.1) we obtain simple 2-persons differential game $DG_{2;T}(\tilde{\mathbf{f}},\mathbf{0},\mathbf{0})$, with linear dynamics:

$$\dot{y}(T,t) = X(T,t)u(t) + Y(T,t)v(t).$$
$$u \in U, v \in V. \tag{A.10}$$
$$J = \|y(T)\|.$$